\documentclass[11pt]{article}

\usepackage{amsmath,amssymb,amsthm}
\usepackage{geometry}
\usepackage{booktabs}
\usepackage{array}
\usepackage[hypertexnames=false]{hyperref}
\usepackage{float}
\usepackage{placeins}
\usepackage[utf8]{inputenc}
\usepackage[T1]{fontenc}
\usepackage{lmodern}
\usepackage{microtype}
\usepackage{longtable}
\usepackage[titletoc,title]{appendix}
\hypersetup{
  hidelinks,
  pdftitle={Arithmetic Bias in the Distribution of Mersenne Prime Exponents and the Divisor Structure of p-1},
  pdfauthor={Jesús Domínguez},
  pdfkeywords={Mersenne primes, divisor structure, arithmetic bias, Wagstaff heuristic, experimental number theory}
}
\newtheorem{lemma}{Lemma}
\newtheorem*{corollary}{Corollary}

\title{\Large\bfseries
Arithmetic Bias in the Distribution of Mersenne Prime Exponents and the Divisor Structure of  \(p-1\)
}

\author{
Jes\'us Dom\'inguez\\[2pt]
\small \texttt{jdomingue601@alumno.uned.es}
}

\date{November 11, 2026\\[2pt]\small arXiv:2603.08994, replacement v8}

\begin{document}

\maketitle

\begin{abstract}
Wagstaff's classical heuristic describes how the probability that
\(2^p-1\) is prime depends on the size of the prime exponent \(p\).
We investigate whether the divisor structure of $p-1$ is associated
with additional variation among exponents of comparable size.
We summarize this structure by \(S(p)=\log\tau(p-1)/\log\log p\), 
where \(\tau\) is the divisor function, and compare regions 
of higher and lower \(S\) within six exponent intervals.
In a Wagstaff-balanced design covering all \(5{,}761{,}450\) prime exponents
\(13\le p<10^8\), \(32\) of the \(47\) known Mersenne-prime exponents lie
in the higher-\(S\) regions, with a conditional two-sided reference value of
\(0.01154\) under the weighted-subset null. These comparisons are exploratory:
the statistic and conditional analyses were selected post hoc, and the present
data do not establish an association beyond the contribution of powers of two
to the divisor count.

We then construct a modular-response benchmark \(H(p)\) from an exactly
characterized cyclotomic order family. The benchmark averages a declared
response over residue configurations satisfying the selected product-one
constraints. Its variation across the full population is summarized by
\(H(p)\approx\theta\tau(p-1)\), with \(\theta\) determined without using
Mersenne-prime labels. Under a restricted effective-channel continuation
hypothesis, relative variation in this score redistributes Wagstaff weight
while preserving its total mass on the finite calibration population.
The resulting rule allocates \(32.28\) and \(14.72\) cases to the higher-
and lower-\(S\) regions, respectively. Extrapolating the same rule to
\(10^8\le p<10^9\) assigns \(73.88\%\) of model weight to the
higher-\(S\) region.
\end{abstract}

\noindent\textbf{MSC 2020:} 11A41, 11N05, 11Y11, 62G10

\noindent\textbf{Keywords:}
Mersenne primes, divisor structure, arithmetic bias,
Wagstaff heuristic, experimental number theory

\section{Introduction}

Wagstaff's heuristic~\cite{Wagstaff} describes the aggregate scale of
Mersenne primality in terms of the exponent \(p\):
\begin{equation}
\mathbb P(M_p\ \mathrm{prime})
\sim C\frac{\log p}{p},
\qquad
M_p=2^p-1.
\label{eq:wagstaff-scale}
\end{equation}
Here \(C\) is Wagstaff's positive leading constant.
The Bateman--Horn heuristic~\cite{BatemanHorn} provides broader context for
prime-value heuristics, although \(2^p-1\) is not a polynomial family and no
direct specialization is asserted here. Primality probabilities denote heuristic
densities over candidate families; a fixed integer is prime or composite.
We use \(b(p)=\log p/p\) as a marginal reference weight and ask whether
arithmetic information in \(p-1\) resolves conditional heterogeneity within
that scale. General background on Mersenne numbers
and their primality testing may be found in Crandall and
Pomerance~\cite{CrandallPomerance}.

We examine the arithmetic structure of \(p-1\) through the parameter
\begin{equation}
S(p)=\frac{\log\tau(p-1)}{\log\log p},
\label{eq:S-intro}
\end{equation}
where \(\tau\) is the divisor function. The motivation is cyclotomic: the
divisors of \(p-1\) index the factors
\[
2^{p-1}-1=\prod_{d\mid p-1}\Phi_d(2).
\]
Here \(\Phi_d\) is the \(d\)-th cyclotomic polynomial. Prime divisors of
the values \(\Phi_d(2)\) select auxiliary moduli for product congruences that every actual
factorization already satisfies. We use this selection to define a
conditional benchmark. Assigning residue one to every factor satisfies
the selected congruences for any formal number of factors; this motivates
the family selection. The congruences supply no new primality exclusion,
and a separately declared response determines the benchmark amplitude.

Related average-order questions belong to the Titchmarsh divisor problem; see,
for example, Drappeau~\cite{Drappeau}. Kova\v{c} and Luca~\cite{KovacLuca}
study the distinct problem of \(\tau(2^k-1)\). Neither line of work addresses
the finite-scale distribution of \(\tau(p-1)\) conditional on \(M_p\) being
prime.

The empirical signal is visible in two complementary designs. Local
comparisons with nearby prime exponents show higher \(S\) values,
supported by rank-based and stratified permutation procedures. These local
analyses use all \(48\) Mersenne-prime exponents \(p\ge13\) in the dataset.

The population-level design uses all \(5{,}761{,}450\) prime
exponents
\[
13\le p<10^8.
\]
Within six disjoint exponent intervals, the population is divided into higher-
and lower-\(S\) regions carrying approximately equal Wagstaff mass. Of the
\(47\) known Mersenne-prime exponents in this range, \(32\) lie in the
higher-\(S\) regions and \(15\) in the lower-\(S\) regions. The primary exact
conditional randomization gives
\[
p_{\mathrm{two}}=0.01154.
\]
This empirical test does not use the structural benchmark developed below.
A post hoc comparison conditioning on \(v_2(p-1)\), the exponent of
\(2\) in \(p-1\), finds no separately detectable odd-part association
(Appendix~\ref{app:order-control}).

We next enumerate the cyclotomic order family generated by \(p-1\).
A declared response to repeated factor residues gives the
\emph{order-selected structural score} \(H(p)\), summarized by the
population approximation \(H(p)\approx\theta\tau(p-1)\).
No Mersenne-primality labels enter these calculations; we call them
\emph{label-independent}.

Controls distinguish the effects of modulus selection, product targets,
and sector weighting on the benchmark.

The averaged response defines a formal density-equivalent index.
We hypothesize that its relative variation, after compression and anchoring,
can redistribute Wagstaff weight. Explicit normalization preserves total
mass on the finite calibration population. Defining the selected index
requires no neutrality assumption for complementary arithmetic, but its
adequacy as a primality model remains unproved.
Section~\ref{sec:reference-anchor} states this continuation hypothesis.
The supplement documents benchmark anatomy, conditional information
limits, the post hoc even-order and \(2\)-adic controls, and
residual-proxy sensitivity.
Additional development experiments remain available in the reproducibility
archive.

The principal compression and \(S=1\) reference give \(32.28/14.72\),
against \(32/15\) observed. The unchanged rule is extrapolated
to \(10^8\le p<10^9\) as an exploratory model implication.

The main contributions are a population-level conditional analysis of
the historical association, an exact cyclotomic order-family framework,
a reproducible structural benchmark, and an explicitly delimited
effective-channel continuation hypothesis.

Sections~2--3 present the arithmetic framework and historical association.
Sections~4--6 define the selected benchmark and its proposed transmission;
Section~7 reports the larger-interval extrapolation. Section~8 and the Appendix
distinguish established results, limitations, and additional audits.

\section{Algebraic constraints and structural interpretation}
\label{sec:algebraic}

\subsection{Algebraic identity}

Let \(M_p=2^p-1\), with \(p\) an odd prime. If \(q\mid M_p\)
is prime, then the multiplicative order of \(2\) modulo \(q\)
equals \(p\). Since this order divides \(q-1\) and \(q\) is odd,
it follows that
\[
q\equiv 1\pmod{2p}.
\]
Thus, we may write
\[
q_i=2pk_i+1,
\qquad
k_i\in\mathbb Z_{\ge 1}.
\]
If \(M_p\) is composite, write
\[
M_p=\prod_{i=1}^{n}q_i,
\qquad n\ge2,
\]
with prime factors counted with multiplicity. Expansion in elementary
symmetric polynomials gives
\[
M_p
=
\prod_{i=1}^{n}(2pk_i+1)
=
1+2p\,A(k_1,\ldots,k_n),
\]
where
\[
A(k_1,\ldots,k_n)
=
\sum_{m=1}^{n}(2p)^{m-1}
e_m(k_1,\ldots,k_n),
\]
and \(e_m\) denotes the \(m\)-th elementary symmetric polynomial.

Since \(M_p=2^p-1\),
\[
A(k_1,\ldots,k_n)
=
\frac{2^{p-1}-1}{p}
=:A_p.
\]
Hence every factorization of \(M_p\) determines a parameter vector
\[
(k_1,\ldots,k_n)
\]
satisfying the exact algebraic relation
\[
A(k_1,\ldots,k_n)=A_p.
\]

\subsection{Cyclotomic layers and modular constraints}

The quantity \(A_p\) inherits the cyclotomic decomposition
\[
2^{p-1}-1
=
\prod_{\substack{d\mid(p-1)\\ d>1}}\Phi_d(2),
\]
where \(\Phi_d(x)\) denotes the \(d\)-th cyclotomic polynomial. Hence the
divisors \(d\mid p-1\) index cyclotomic layers whose prime divisors provide
auxiliary moduli associated with the divisor structure of \(p-1\).

Let \(a\ne p\) be a prime divisor of \(A_p\), and put
\[
d=\operatorname{ord}_a(2).
\]
Then
\[
d>1,
\qquad
d\mid p-1,
\]
so \(a\) belongs to the cyclotomic order layer indexed by \(d\). Moreover,
\[
a\mid A_p
\quad\Longrightarrow\quad
2^{p-1}\equiv1\pmod a,
\]
and therefore
\[
M_p
=
2(2^{p-1}-1)+1
\equiv1\pmod a.
\]
Thus every factorization
\[
M_p=\prod_{i=1}^{n}q_i
\]
satisfies the exact joint residue constraint
\begin{equation}
\prod_{i=1}^{n}q_i\equiv1\pmod a.
\label{eq:joint-residue-constraint}
\end{equation}
Equivalently, using \(q_i=2pk_i+1\),
\[
\prod_{i=1}^{n}(2pk_i+1)\equiv1\pmod a.
\]

The joint congruence is a necessary consequence of the product identity
for every realized factorization, not an independent probabilistic exclusion.
Its role here is to select the product-one residue fibers associated with
the divisor structure of \(p-1\) for the benchmark of
Section~\ref{sec:structural-density-model}.
No primality modulation follows from this identity alone.

\subsection{Cyclotomic complexity}

Write \(\omega(n)\) for the number of distinct prime divisors of \(n\).
The divisors \(d\mid p-1\) index \(\tau(p-1)-1\) nontrivial
cyclotomic layers. Thus \(\tau(p-1)\), unlike \(\omega(p-1)\),
retains the contribution of prime-power multiplicities.
Zsigmondy's theorem turns this layer count into a rigorous lower
bound on the number of distinct prime contributions to \(A_p\).

\begin{lemma}[Cyclotomic-layer lower bound]
For every odd prime \(p\),
\[
\omega(A_p)\geq \tau(p-1)-3.
\]
\end{lemma}

\begin{proof}
A closely related general Bang--Zsigmondy bound,
\(\omega(2^k-1)\geq\tau(k)-2\), is used by Kova\v{c} and
Luca~\cite{KovacLuca}. Because the present statement concerns the quotient
\(A_p=(2^{p-1}-1)/p\), we record the short direct argument. Zsigmondy's
theorem implies that, for each divisor
\(d\mid p-1\) with \(d>1\) and \(d\ne6\), the number \(2^d-1\)
has a primitive prime divisor~\cite{Zsigmondy}.  A primitive divisor
belonging to the layer \(d\) cannot be primitive for a different layer,
so the eligible divisors \(d\) supply distinct prime factors of
\(2^{p-1}-1\).  The index \(d=1\) gives \(\Phi_1(2)=1\) and therefore
contributes no prime divisor, the
exceptional layer \(d=6\) may be absent, and division by \(p\) in
\(A_p=(2^{p-1}-1)/p\) can remove at most one further distinct prime.
At most three elements are therefore lost from the \(\tau(p-1)\)
cyclotomic indices.  If \(p\) is Wieferich to base \(2\), division by
\(p\) need not remove it, which can only improve the bound.
\end{proof}

Thus, up to an additive constant, \(\tau(p-1)\) gives a lower bound on
the number of distinct cyclotomic prime contributions to \(A_p\). This bound
concerns the full prime-divisor set of \(A_p\). The effective census retains
only \(a<p\), so its finite-scale scaling is examined empirically in
Appendix~\ref{app:empirical}.

\subsection{Normalized divisor structure}

Recall the normalized divisor statistic
\[
S(p)=\frac{\log \tau(p-1)}{\log\log p}.
\]

The normalization by \(\log\log p\) keeps the statistic on a comparable
scale across the tested exponent range. In particular, the reference value \(S=1\)
corresponds algebraically to
\[
\tau(p-1)=\log p.
\]
This reference scale will later define the reduced model's relative
factor; it is not assumed to represent the mean or
typical value of \(S\) among prime exponents.

The choice of \(\tau(p-1)\) is suggested directly by the cyclotomic
decomposition, since prime-power multiplicities generate additional
divisor layers not recorded by \(\omega(p-1)\). Writing
\[
p-1=\prod_{i=1}^{\omega(p-1)}\ell_i^{e_i},
\]
we have
\[
\log\tau(p-1)
=
\sum_{i=1}^{\omega(p-1)}\log(e_i+1).
\]
Classical probabilistic number theory places \(\omega(n)\) on the
\(\log\log n\) scale, with the Erd\H{o}s--Kac theorem describing its
fluctuations about that scale. Since the exponents \(e_i\) are
typically small, the sum above suggests comparison of \(\log\tau(p-1)\)
with \(\log\log p\) at the finite scales considered here~\cite{HardyWright,Tenenbaum,ErdosKac}. Appendix~A compares \(S\)
with normalized statistics based separately on distinct-prime count
and multiplicity.

\medskip

\subsection{CRT residue assignments and singleton coordinates}
\label{sec:residue-solutions}

Let \(a\ne p\) be a prime divisor of \(A_p\). Since
\(A(k_1,\ldots,k_n)=A_p\),
\[
A(k_1,\ldots,k_n)\equiv0\pmod a.
\]
After setting
\[
u_i\equiv2pk_i+1\pmod a,
\]
this becomes
\[
u_1\cdots u_n\equiv1\pmod a,
\qquad u_i\in\mathbb F_a^\times.
\]
The ordered product-one reference space, or fiber, is therefore
\begin{equation}
\Omega_{a,n}
=
\left\{(u_1,\ldots,u_n)\in(\mathbb F_a^\times)^n:
u_1\cdots u_n=1\right\},
\qquad
|\Omega_{a,n}|=(a-1)^{n-1}.
\label{eq:omega-an}
\end{equation}
The coordinate-wise CRT calculation uses ordered residue tuples. A sector
has a fixed number \(n\) of prime factors; its geometric weight is defined
separately in Section~\ref{sec:sector-measure}.

Every candidate prime factor already lies in the class \(1\pmod{2p}\). For
each \(u\in\mathbb F_a^\times\), the Chinese remainder theorem gives a unique
reduced class \(r_u\pmod{2pa}\) satisfying
\[
r_u\equiv1\pmod{2p},
\qquad
r_u\equiv u\pmod a.
\]
We call a specified auxiliary residue class a \emph{channel}.
The prime number theorem for arithmetic progressions gives~\cite{Tenenbaum}
\begin{equation}
\lim_{x\to\infty}
\frac{\pi(x;2pa,r_u)}{\pi(x;2p,1)}
=
\frac{\varphi(2p)}{\varphi(2pa)}
=
\frac1{a-1}.
\label{eq:relative-channel-density}
\end{equation}
Here \(\pi(x;m,r)\) counts primes at most \(x\) congruent to \(r\) modulo
\(m\), and \(\varphi\) is Euler's totient function. Thus \(1/(a-1)\)
is the limiting relative density of a fixed \(u\)-channel inside the
classical admissible-prime support. A normalized weighted average of
these equal class-density ratios has the same value; it is not the
density of a union of channels.

For \(\mathbf u\in\Omega_{a,n}\), call coordinate \(u_i\) a \emph{singleton}
when no other coordinate has the same value. The probability below concerns one specified coordinate, rather than
the absence of collisions in the whole assignment.

\begin{lemma}[Exact singleton-coordinate probability]
\label{lem:singleton-probability}
Let \(a\) be an odd prime, \(n\ge2\), and put \(m=a-1\). Under the uniform
measure on \(\Omega_{a,n}\), the probability that a specified coordinate is a
singleton is
\begin{equation}
P_{\mathrm{sing}}(a,n)
=
\frac{(m-1)^{n-1}
+(-1)^{n-1}\bigl(\gcd(n,m)-1\bigr)}{m^{n-1}}.
\label{eq:psing}
\end{equation}
\end{lemma}

\begin{proof}
By symmetry it is enough to consider the first coordinate. Put
\[
G:=\mathbb F_a^\times,
\qquad |G|=m,
\]
and fix \(u_1=x\in G\). Write
\[
u_j=xy_j,
\qquad 2\le j\le n.
\]
Then \(u_1\) is a singleton exactly when \(y_j\ne1\) for all \(j\ge2\),
while the constraint
\[
u_1\cdots u_n=1
\]
becomes
\[
y_2\cdots y_n=x^{-n}.
\]

Put \(h=n-1\), and for \(g\in G\) define
\[
N_h(g)
=
\#\left\{
(z_1,\ldots,z_h)\in(G\setminus\{1\})^h:
z_1\cdots z_h=g
\right\}.
\]
Thus, for fixed \(u_1=x\), the number of admissible assignments for which
\(u_1\) is a singleton is \(N_h(x^{-n})\).

By character orthogonality,
\[
N_h(g)
=
\frac1m
\sum_{\chi\in\widehat G}
\chi(g^{-1})
\left(
\sum_{z\in G\setminus\{1\}}\chi(z)
\right)^h.
\]
For the trivial character the inner sum is \(m-1\), whereas for every
nontrivial character it is \(-1\). Hence
\[
N_h(g)
=
\frac1m
\left[
(m-1)^h
+
(-1)^h
\sum_{\substack{\chi\in\widehat G\\ \chi\ne1}}
\chi(g^{-1})
\right].
\]

Summing over \(x\in G\), with \(g=x^{-n}\), gives
\[
\sum_{x\in G}N_h(x^{-n})
=
(m-1)^h
+
\frac{(-1)^h}{m}
\sum_{\substack{\chi\in\widehat G\\ \chi\ne1}}
\sum_{x\in G}\chi(x)^n.
\]
Now
\[
\sum_{x\in G}\chi(x)^n
=
\begin{cases}
m, & \chi^n=1,\\
0, & \chi^n\ne1.
\end{cases}
\]
Because \(\widehat G\) is cyclic of order \(m\), the equation
\(\chi^n=1\) has exactly \(\gcd(n,m)\) solutions, one of which is the
trivial character. Therefore the number of assignments for which \(u_1\) is
a singleton is
\[
(m-1)^h+(-1)^h\bigl(\gcd(n,m)-1\bigr).
\]
Since \(h=n-1\) and \(|\Omega_{a,n}|=m^{\,n-1}\),
\[
P_{\mathrm{sing}}(a,n)
=
\frac{
(m-1)^{n-1}
+
(-1)^{n-1}\bigl(\gcd(n,m)-1\bigr)
}{m^{\,n-1}},
\]
as claimed.
\end{proof}

For \(n=2\), formula~\eqref{eq:psing} reduces to
\[
P_{\mathrm{sing}}(a,2)=\frac{a-3}{a-1}.
\]
Indeed, the ordered assignments are \((u,u^{-1})\); only \(u=\pm1\) are
diagonal. The same formula treats all \(n\ge2\) without singling out the
two-factor sector.

For later asymptotics, when \(a\gg n\), equation~\eqref{eq:psing} gives
\begin{equation}
P_{\mathrm{sing}}(a,n)
=
1-
\frac{n-1+\mathbf 1_{\{n=2\}}}{a-1}
+O\!\left(\frac{n^2}{a^2}\right).
\label{eq:psing-large-a}
\end{equation}
The exceptional first-order gcd correction at \(n=2\) is displayed explicitly;
for \(n\ge3\), the remaining gcd term is absorbed by the stated error. Hence
the probability that a specified coordinate is repeated becomes small for moduli large relative to the sector length.

\section{Historical association}
\label{sec:empirical}

Appendix~A compares \(S(p)\) for known Mersenne-prime exponents with
nearby prime controls of comparable size, using local percentiles,
stratified conditional logistic regression, and stratified permutation tests.

In the pooled local comparisons, known Mersenne-prime exponents have
elevated \(S(p)\) relative to nearby prime controls. With \(W=5000\)
controls on each side, truncated at the lower boundary, the mean percentile
rank is \(\overline{\pi}_{5000}=0.648\), above the \(0.5\) center of the
uniform within-stratum reference.

For this window, uniform within-stratum randomization gives the approximate
reference value \(p_{\mathrm{perm}}=0.0012\). Overlap between choice sets
and the wider ranges of the earliest boundary-truncated sets make these
local values descriptive. The population analysis below avoids local
windows and provides the primary comparison under the specified
Wagstaff-weighted conditional null.

Within each of six disjoint exponent intervals, condition on the
observed case count and assign each subset probability proportional to
the product of its Wagstaff weights. This gives
\[
X_+=32\ \text{of}\ 47,
\qquad
p_{\mathrm{two}}=0.01154
\]
under the exact conditional weighted-subset null.  The two-sided value
is the primary reference value. Directional and
congruence-sensitive variants are reported with the complete
population-level calculation in Appendix~\ref{app:wagstaff_balanced}.

For \(W=5000\), the pooled local-control analysis gives a standardized
mean difference \(d\approx0.56\), using the average-variance convention
specified in Appendix~\ref{app:local-summary}.

In these historical comparisons, the recorded Mersenne-prime exponents
have richer divisor structure in \(p-1\) than the nearby prime controls.

Appendices~\ref{app:exponent-data}--\ref{app:local-robustness} detail the
local analyses; Appendix~\ref{app:wagstaff_balanced} gives the primary
population test. The post hoc \(v_2(p-1)\)-conditioned control in
Appendix~\ref{app:order-control} is discussed in Section~\ref{sec:limitations}.

The following sections construct a structural score from the modular
constraints generated by \(p-1\) and state the additional hypothesis needed
to map that score onto Wagstaff-scale allocations.

\medskip

\section{Wagstaff scale and the auxiliary-channel benchmark}
\label{sec:wagstaff-refinement}

\subsection{Classical Wagstaff density}

A compact sieve reconstruction of Wagstaff's leading scale begins with
the random-integer scale
\[
\frac{1}{\log M_p}\sim\frac{1}{p\log2}
\]
and uses the fact that every prime divisor of \(M_p\) belongs to the reduced
class \(1\pmod{2p}\). Relative to an unrestricted integer, the exclusion of
inadmissible small prime divisors supplies the Mertens correction
\[
\prod_{\ell<2p}\left(1-\frac1\ell\right)^{-1}
\sim e^\gamma\log(2p).
\]
Here \(\gamma\) is Euler's constant. This gives the classical heuristic scale
\[
\mathbb P(M_p\ \mathrm{prime})
\approx
\frac{e^\gamma\log(2p)}{p\log2}.
\]
The factor restriction \(q\equiv\pm1\pmod8\) distinguishes the two
cases of \(p\bmod4\), giving Wagstaff's more precise numerator
\[
\mathbb P(M_p\ \mathrm{prime})
\approx
\frac{e^\gamma\log(a(p)p)}{p\log2},
\qquad
a(p)=
\begin{cases}
2,&p\equiv3\pmod4,\\
6,&p\equiv1\pmod4.
\end{cases}
\]
Here \(a(p)\in\{2,6\}\) is Wagstaff's size coefficient, distinct from the
auxiliary prime modulus \(a\) used below.
Wagstaff's original argument~\cite[pp.~385--388]{Wagstaff} is more specific
than this compact sieve reconstruction. For fixed \(k\), it combines a
Bateman--Horn estimate for simultaneous primality of \(p\) and
\(q=2kp+1\) with a heuristic divisibility probability of approximately
\(1/k\) when \(q\) is prime and \(k\equiv0\) or \(-p\pmod4\),
followed by a Poisson zero-count argument and a correction to match
prime-number-theorem normalization. The reconstruction above retains the
leading logarithmic scale. The formal auxiliary-channel index below does
not inherit this divisibility law; its transmission to primality remains
a separate heuristic hypothesis.

The leading form is~\eqref{eq:wagstaff-scale}, with \(C=e^\gamma/\log2\).
We use
\(b(p)=\log p/p\) as the primary marginal reference weight and retain
\(\log(a(p)p)/p\) as a sensitivity check. No structural conditioning is introduced
at this stage.

The benchmark below describes selected auxiliary arithmetic. Its
proposed primality interpretation is stated in
Section~\ref{sec:reference-anchor}.

\subsection{Fixed-channel densities and effective-density motivation}
\label{sec:structural-density-bridge}

The classical support is factorwise: every candidate prime factor lies
in \(q\equiv1\pmod{2p}\). The auxiliary condition
\(\prod_iq_i\equiv1\pmod a\), for an active prime modulus \(a<p\) dividing
\(A_p\), instead constrains
the joint residue assignment. For a fixed \(u\in\mathbb F_a^\times\),
CRT combines \(q\equiv1\pmod{2p}\) and \(q\equiv u\pmod a\) into
one reduced class \(r_u\pmod{2pa}\). Writing \(\delta(r\bmod m)=1/\varphi(m)\) for the limiting
proportion of primes in a reduced residue class,
Equation~\eqref{eq:relative-channel-density} gives
\begin{equation}
\frac{\delta(r_u\bmod 2pa)}{\delta(1\bmod 2p)}
=\frac{\varphi(2p)}{\varphi(2pa)}
=\frac1{a-1}.
\label{eq:bridge-channel-density}
\end{equation}
This fixed-modulus limiting ratio is independent of \(u\): the identity
channel has no density advantage. Averaging the equal class-density
ratios preserves their value, whereas a union of distinct channels adds
their densities.

Repeated coordinates impose a shared-channel condition: if
\(u_i=u_j\), the corresponding candidate factors satisfy
\(q_i\equiv q_j\pmod{2pa}\). This constrains the joint factor
configuration without requiring the shared residue to be one.
The auxiliary modulus labels these channels; it is not a candidate
divisor of \(M_p\). Section~\ref{sec:local-survival} specifies how
shared-channel membership enters the declared response.

\paragraph{Conditional strengthening of factor support.}
For an active prime modulus $a<p$, the product-one CRT
fiber contains the all-identity assignment. In this
particular configuration,
\[
q_i\equiv1\pmod{2p},
\qquad
q_i\equiv1\pmod a,
\]
and hence, by CRT,
\[
q_i\equiv1\pmod{2pa},
\qquad
q_i\ge 2pa+1.
\]

This provides a conditional strengthening of the
classical admissible-factor support. In Wagstaff's
leading logarithmic heuristic, a hypothetical
extension of the factor-size threshold from order
$p$ to order $ap$ suggests the formal analogy
\[
\frac{C\log p}{p}
\quad\rightsquigarrow\quad
\frac{C\log(ap)}{p},
\]
where $a$ is the auxiliary prime modulus, distinct
from Wagstaff's size coefficient $a(p)$.

Unlike Wagstaff's universal admissible-factor restrictions,
the stronger auxiliary congruence holds only within the
all-identity assignment. The product-one constraint does not
impose this assignment. The logarithmic analogy is therefore
motivational rather than a derived sieve correction: it
establishes neither a global factor-size cutoff nor a new
primality law. The use of the same leading constant $C$ is
illustrative rather than a derived normalization.

\paragraph{Effective-density averaging.}
The all-identity assignment is only one of the
residue configurations compatible with the
product-one constraint. To represent a response
incorporating different configurations, we adopt
an effective-density averaging convention.

For reference components with weights
$\mu_j\ge0$, $\sum_j\mu_j=1$, and declared
relative responses $0<\rho_j\le1$, define
\begin{equation}
\rho_{\mathrm{eff}}
=
\sum_j\mu_j\rho_j.
\label{eq:density-mixture-principle}
\end{equation}
The responses are averaged before any
inverse-density or logarithmic transformation.

This operation must be distinguished from
averaging the genuine Dirichlet densities of
residue classes. Equation~\eqref{eq:bridge-channel-density} assigns the same
limiting relative density $1/(a-1)$ to each
fixed auxiliary channel; averaging those equal
densities cannot by itself produce an additional
restriction.

A nontrivial effective response therefore requires
a separately specified model of residue
configurations. Section~\ref{sec:local-survival}
introduces this response, and
Section~\ref{sec:effective-aggregation}
specifies its averaging and aggregation.
Its algebraic interpretation is developed in
Section~\ref{sec:identity-boundary}, before
the proposed continuation to Wagstaff's
primality scale.

\section{Order-family benchmark: local response, sector measure, and aggregation}
\label{sec:structural-density-model}

The structural benchmark is evaluated for prime exponents \(p\ge13\),
for which the sector range defined below is nonempty. The reference
sectors are formal and are evaluated even when \(M_p\) is prime.
Section~\ref{sec:algebraic} gives the algebraic activation rule and the residue spaces
\(\Omega_{a,n}\); Section~\ref{sec:structural-density-bridge} supplies the
fixed-channel density ratio \(1/(a-1)\). We now specify the factorwise
response and the sector and modulus aggregation rules. The relative factor
transmitted heuristically to the Wagstaff scale is defined after
population-level compression in Section~\ref{sec:reference-anchor}.

\subsection{Effective moduli and exact census support}
\label{sec:effective-moduli}

For a prime exponent \(p\), define
\begin{equation}
\mathcal A_p
=
\{a<p:a\ \text{odd prime and }\operatorname{ord}_a(2)\mid p-1\},
\qquad
L(p)=|\mathcal A_p|.
\label{eq:active-set}
\end{equation}
Since \(a<p\), division by \(p\) does not affect divisibility by \(a\), and
the equivalences
\[
a\mid A_p
\quad\Longleftrightarrow\quad
2^{p-1}\equiv1\pmod a
\quad\Longleftrightarrow\quad
\operatorname{ord}_a(2)\mid p-1
\]
are exact.

For a fixed numerical cutoff \(B\), let
\(L_B(p):=\#\{a\in\mathcal A_p:a<B\}\), and define
\[
c_d(X)=\#\{a<X:a\ \text{odd prime},\ \operatorname{ord}_a(2)=d\}.
\]
Then the captured count satisfies the exact identity
\begin{equation}
L_B(p)
=
\sum_{d\mid p-1}c_d\!\bigl(\min(B,p)\bigr).
\label{eq:order-capture}
\end{equation}
For \(p\ge B\), this reduces to \(L_B(p)=\sum_{d\mid p-1}c_d(B)\).
Equation~\eqref{eq:order-capture} links the divisors of \(p-1\) to the count
of active prime moduli grouped by multiplicative order. It also gives
an analogy with the multiplicative-order criterion in Pollard's
\(p-1\) method, but it does not by itself imply a primality correction.

Appendix~\ref{app:empirical} uses this identity to audit the active-modulus
census. The following subsection examines its empirical count scaling.

\subsection{Effective-modulus count scaling}
\label{sec:count-scaling}

Equation~\eqref{eq:order-capture} suggests a positive relation between
the active-modulus count and the divisor function, without implying
linearity. The complete census of \(5{,}761{,}450\) prime exponents measures
the accuracy and scale dependence of a through-origin approximation for
\(L(p)\), with \(a<p\) as in Equation~\eqref{eq:active-set}.
For an exponent interval \(I\), write
\(\mathcal P_I=\{p\in I:p\text{ prime}\}\); all interval sums run over
these prime exponents. We fit the through-origin approximation
\begin{equation}
L(p)\simeq \kappa_I\tau(p-1)
=\kappa_I(\log p)^{S(p)}.
\label{eq:L-regression-new}
\end{equation}
Table~\ref{tab:count-scaling} reports the fitted slopes. The Pearson and
Spearman correlations
remain above \(.958\) and \(.961\), respectively, in every decade.

\begingroup
\setlength{\intextsep}{3pt}
\setlength{\abovecaptionskip}{4pt}
\begin{table}[H]
\centering
\small
\begin{tabular}{lrrrr}
\toprule
Exponent interval & Population & \(\kappa_I\) & Pearson \(r\) & Spearman \(\rho\) \\
\midrule
\([13,10^3)\) & 163 & .523328 & .9589 & .9613 \\
\([10^3,10^4)\) & 1,061 & .582290 & .9671 & .9622 \\
\([10^4,10^5)\) & 8,363 & .621716 & .9728 & .9652 \\
\([10^5,10^6)\) & 68,906 & .689001 & .9803 & .9711 \\
\([10^6,10^7)\) & 586,081 & .719268 & .9834 & .9732 \\
\([10^7,10^8)\) & 5,096,876 & .749332 & .9857 & .9749 \\
\midrule
Global & 5,761,450 & .746752 & .9855 & .9747 \\
\bottomrule
\end{tabular}
\caption{Population regression \(L\simeq \kappa_I\tau(p-1)\) from the complete
active-modulus census. The slopes and correlations are descriptive
summaries; no Mersenne-prime label is used.}
\label{tab:count-scaling}
\end{table}
\endgroup

Equation~\eqref{eq:L-regression-new} approximates the number of active
moduli before their heterogeneous local responses are included.

\subsection{Maximum sector length and unordered sector measure}
\label{sec:sector-measure}

Every prime divisor \(q\mid M_p\) satisfies \(q\equiv1\pmod{2p}\).
Since the order of \(2\) modulo \(q\) is odd, \(2\) is a quadratic residue
modulo \(q\), so \(q\equiv\pm1\pmod8\). Writing \(q=1+2pk\), the least
admissible positive \(k\) is \(1\) for \(p\equiv3\pmod4\) and \(3\) for
\(p\equiv1\pmod4\). These congruences give the lower bound
\begin{equation}
Q_{\min}(p)=
\begin{cases}
2p+1,&p\equiv3\pmod4,\\
6p+1,&p\equiv1\pmod4.
\end{cases}
\label{eq:qmin}
\end{equation}
The bound does not require \(Q_{\min}(p)\) to be prime. Define
\begin{equation}
R_p=\frac{\log(2^p-1)}{\log Q_{\min}(p)},
\qquad
N_{\max}(p)=\lfloor R_p\rfloor.
\label{eq:Rp-Nmax}
\end{equation}
A composite factorization into \(n\) prime factors satisfies
\(2\le n<R_p\), with
\[
N_{\max}(p)\asymp\frac{p}{\log p}.
\]

We assign a geometric weight to each possible sector \(n\).
For a factorization \(M_p=q_1\cdots q_n\), introduce the normalized
log-size variables
\[
x_i=\frac{\log(q_i/Q_{\min})}{\log Q_{\min}}\ge0.
\]
Their sum is fixed:
\[
\sum_{i=1}^n x_i=R_p-n.
\]
We equip this simplex with the coordinate measure
\(dx_1\cdots dx_{n-1}\), obtained after eliminating \(x_n\).
This is not the Euclidean surface measure on the embedded
hyperplane. With this convention, the corresponding simplex
has volume
\[
V_n(p)=\frac{(R_p-n)^{n-1}}{(n-1)!}.
\]

Permuting the \(n\) factors does not change the factorization.
Quotienting generic ordered points by permutations contributes
\(1/n!\); coincident-coordinate strata have zero measure under
this continuous measure. We therefore define
\begin{equation}
\widetilde w_n(p)
=
\frac{(R_p-n)^{n-1}}{(n-1)!\,n!},
\qquad
w_n(p)=
\frac{\widetilde w_n(p)}
{\displaystyle\sum_{2\le m<R_p}\widetilde w_m(p)}.
\label{eq:wn}
\end{equation}

These weights define the benchmark's geometric allocation across
factor-number sectors. They are a modeling choice, not an asserted
probability distribution for the actual number of prime factors
of \(M_p\). The factor \(1/n!\) applies only to this sector
allocation; the within-sector CRT calculations continue to use
the ordered space \(\Omega_{a,n}\). Appendix~\ref{app:sector-prior-sensitivity}
examines the sensitivity to the sector weights.

For numerical evaluation, we compute
\[
\ell_n
=(n-1)\log(R_p-n)-\log\Gamma(n)-\log\Gamma(n+1),
\]
and normalize using log-sum-exp. The mode of the weights is of
order \(\sqrt{R_p}\), well below \(N_{\max}(p)\).

\subsection{Finite-range singleton diagnostic}
\label{sec:singleton-diagnostic}

The exact geometry can be compared with realized factorizations without using
the comparison to estimate any model coefficient. For a completely factored
composite
\[
M_p=\prod_{i=1}^nq_i,
\qquad n\ge3,
\]
and each active \(a<p\), put \(u_i=q_i\bmod a\) and let
\[
N_{\mathrm{sing}}^{\mathrm{obs}}(p,a)
=\#\{i:u_i\ne u_j\ \text{for all }j\ne i\}.
\]
Under the uniform measure on \(\Omega_{a,n}\),
Lemma~\ref{lem:singleton-probability} gives the exact mean
\begin{equation}
\mathbb E_\Omega N_{\mathrm{sing}}(a,n)
=nP_{\mathrm{sing}}(a,n).
\label{eq:expected-singletons}
\end{equation}

The diagnostic uses all \(143\) completely factored composite \(M_p\) with
prime \(p<1200\) and \(n\ge3\), producing \(833\) pairs \((p,a)\). The distribution of \(N_{\mathrm{sing}}\) for each pair is counted exactly
from occupancy partitions, retaining the dependence among singleton
indicators within one residue vector. Within each exponent, CRT factors the uniform
reference over distinct active moduli. For the pooled statistic, we take
the reference law to be the product of these exponent-level laws and
convolve the pair distributions. The fixed-sector product space and measure
are stated in Section~\ref{sec:effective-aggregation}.

\begin{table}[H]
\centering
\small
\begin{tabular}{lrrrrr}
\toprule
Exponent range & \(\#p\) & \(\#(p,a)\) & Observed & Expected & \(z\) \\
\midrule
\(p<500\) & 53 & 243 & 679 & 675.577 & .260 \\
\(500\le p<1000\) & 64 & 406 & 1,269 & 1,264.245 & .305 \\
\(1000\le p<1200\) & 26 & 184 & 621 & 616.611 & .415 \\
\midrule
\(p<1200\) & 143 & 833 & 2,569 & 2,556.433 & .547 \\
\bottomrule
\end{tabular}
\caption{Observed and exact-reference totals of singleton coordinates. For the
pooled statistic, the centered two-sided product-null probability is
\(0.5869\); the three disjoint-range values are \(0.8199\), \(0.7732\), and
\(0.7055\).}
\label{tab:singleton-diagnostic}
\end{table}

The pooled singleton total shows no detectable departure from the
uniform CRT reference in this selected sample of small, completely factored
composites. This diagnostic concerns the singleton total, not the full
occupancy geometry. It tests neither the declared response assignment
nor its subsequent transmission to primality.
The fixed-channel density \(1/(a-1)\) is supplied independently by the
limiting prime-density ratio in arithmetic progressions; cross-sector
weighting, modular aggregation, and the Wagstaff transmission enter at
later stages of the construction.

Secondary deficit and collision summaries likewise show no detectable
departure from the CRT reference and are reported in
Appendix~\ref{app:singleton-audit}. They reuse the same factor sample and
are supporting occupancy diagnostics, not independent tests.

\subsection{Declared factorwise response at fixed sector}
\label{sec:local-survival}

For each active modulus \(a\) and fixed sector \(n\), give every factor
coordinate equal reference weight \(1/n\). For
\[
\mathbf u=(u_1,\ldots,u_n)\in\Omega_{a,n},
\]
let
\[
R(\mathbf u)
:=
\#\{j:\exists\,k\neq j\text{ with }u_j=u_k\}.
\]

We summarize the joint residue geometry through the fraction of
coordinates that repeat. This fraction, \(R(\mathbf u)/n\),
receives the fixed-channel factor \(1/(a-1)\) from
Section~\ref{sec:structural-density-bridge}; the singleton fraction receives
the neutral factor \(1\). These assignments specify the response, rather
than follow from the CRT product constraint or the prime-density ratio.
They do not define a joint prime-density law for complete factorizations.
Thus
\[
\rho_{a,n}(\mathbf u)
=\left(1-\frac{R(\mathbf u)}n\right)
+\frac{R(\mathbf u)}{n(a-1)}.
\]
For the all-identity assignment, \(R=n\) and the response is \(1/(a-1)\).
Other assignments can attain this endpoint: repeated groups need not
occupy the identity class or all share the same class. The equivalent
index represents the averaged scalar response, not a replacement of
these configurations by an all-identity vector.

Under the uniform reference measure on \(\Omega_{a,n}\), exchangeability gives
\[
\mathbb E_{\Omega_{a,n}}
\left[
\frac{R}{n}
\right]
=
1-P_{\rm sing}(a,n),
\]
because \(P_{\rm sing}(a,n)\) is the probability that a uniformly sampled
factor coordinate belongs to a singleton residue class. Averaging the
retained relative response therefore gives
\begin{equation}
\begin{aligned}
s_{a,n}
&=
\mathbb E_{\Omega_{a,n}}
\bigl[\rho_{a,n}(\mathbf u)\bigr]\\
&=
1-
\frac{\mathbb E_{\Omega_{a,n}}R}{n}
+
\frac{\mathbb E_{\Omega_{a,n}}R}{n(a-1)}\\
&=
P_{\mathrm{sing}}(a,n)
+
\frac{1-P_{\mathrm{sing}}(a,n)}{a-1}.
\end{aligned}
\label{eq:survival-an}
\end{equation}

Equation~\eqref{eq:survival-an} is the exact expectation of the specified
response under the uniform reference measure.

Since \(0\le P_{\mathrm{sing}}(a,n)\le1\), convex sector averaging preserves
the bounds
\[
\frac1{a-1}\le s_{a,n}\le1,\qquad
\frac1{a-1}\le\overline s_a(p)\le1,\qquad
1\le F_a(p)\le a-1,
\]
where the sector average \(\overline s_a(p)\) and its inverse \(F_a(p)\) are defined
in Section~\ref{sec:effective-aggregation}. These endpoints bound the declared
response; an averaged factor need not attain the upper endpoint.

For example, \(a=7\), \(n=2\) gives \(P_{\mathrm{sing}}=2/3\),
\(s_{7,2}=13/18\), and sector-level inverse response \(18/13\).
This illustrates why an equivalent index need not be an integer.

\paragraph{Target averaging within the declared response.}
Put \(m=a-1\), and write \(s_{a,n}(t)\) for the same declared response
averaged on the product-target \(t\) fiber, with \(t\in\mathbb F_a^\times\).
Averaging uniformly over targets gives
\(P_{\mathrm{sing}}^{\mathrm{gen}}(a,n)=(1-1/m)^{n-1}\), and hence
\[
s_{a,n}^{\mathrm{gen}}=\frac1m+\left(1-\frac1m\right)^n,
\qquad
s_{a,n}(1)-s_{a,n}^{\mathrm{gen}}
=(-1)^{n-1}\frac{(m-1)[\gcd(n,m)-1]}{m^n}
=O_n(m^{1-n}).
\]
Supplement~S.1 derives these formulas from the
general-target geometry. The correction need not be small for low sectors:
at \(n=2\) it has the same \(1/m\) order as the leading response change.
The population control below measures its contribution after sector
averaging, without inferring target independence solely from the asymptotic.

\subsection{Effective modular aggregation and the boundary
\texorpdfstring{\(a<p\)}{a<p}}
\label{sec:effective-aggregation}

For fixed \(n\), the modular product structure of the reference space is
exact. Put \(Q_p:=\prod_{a\in\mathcal A_p}a\), and define
\[
\Omega_{Q_p,n}
:=
\left\{(U_1,\ldots,U_n)\in
\bigl((\mathbb Z/Q_p\mathbb Z)^\times\bigr)^n:
\prod_{i=1}^n U_i=1\right\}.
\]
The active moduli are distinct primes, so the Chinese remainder theorem gives
\[
(\mathbb Z/Q_p\mathbb Z)^\times
\simeq
\prod_{a\in\mathcal A_p}\mathbb F_a^\times,
\qquad
\Omega_{Q_p,n}
\simeq
\prod_{a\in\mathcal A_p}\Omega_{a,n}.
\]
Here \(Q_p\) is a modulus for the reference space, not a lower bound
on individual factors or an effective cutoff for complete factorizations.
Writing \(\mu_{Q_p,n}\) and \(\mu_{a,n}\) for the uniform measures on these
finite spaces, respectively, the latter identification also gives
\[
\mu_{Q_p,n}
=
\bigotimes_{a\in\mathcal A_p}\mu_{a,n}.
\]
If, in addition, the joint benchmark response at fixed \(n\) is represented
multiplicatively by
\[
\rho_{Q_p,n}\bigl((\mathbf u_a)_{a\in\mathcal A_p}\bigr)
:=
\prod_{a\in\mathcal A_p}\rho_{a,n}(\mathbf u_a),
\]
then its expectation factorizes exactly under this reference measure:
\[
\mathbb E_{\Omega_{Q_p,n}}[\rho_{Q_p,n}]
=
\prod_{a\in\mathcal A_p}
\mathbb E_{\Omega_{a,n}}[\rho_{a,n}]
=
\prod_{a\in\mathcal A_p}s_{a,n}.
\]
CRT establishes the product reference space and uniform measure; the
multiplicative response is a benchmark choice. There is a further choice
in the order of averaging: the principal benchmark averages within each
modulus before aggregation, whereas the shared-sector alternative retains
a common global \(n\). Appendix~\ref{app:allocation-audit} compares these
constructions and finds nearly unchanged conditioned allocations.
For the principal benchmark, the normalized sector weights give
\begin{align}
\overline s_a(p)
&=
\sum_{2\leq n<R_p}w_n(p)s_{a,n},
\label{eq:survival-average}\\
F_a(p)
&=
\overline s_a(p)^{-1}.
\label{eq:Fa-correct}
\end{align}
Equivalently,
\begin{equation}
F_a(p)
=
\left[
\frac{1}{a-1}
+
\left(
1-\frac{1}{a-1}
\right)
\sum_n w_n(p)P_{\rm sing}(a,n)
\right]^{-1}.
\label{eq:Fa-explicit}
\end{equation}

At the modular level, \(\overline s_a(p)\) is the sector-averaged
factorwise response associated with modulus \(a\). Applying the declared
modular product to these separately marginalized responses, define
\[
\rho_{\mathrm{cyc,cens}}(p)
:=
\prod_{a\in\mathcal A_p}\overline s_a(p).
\]
The corresponding inverse-response index is
\begin{equation}
\lambda_{\mathrm{cens}}(p)
=
\rho_{\mathrm{cyc,cens}}(p)^{-1}
=
\prod_{a\in\mathcal A_p}F_a(p),
\qquad
H(p)
=
\log\lambda_{\mathrm{cens}}(p)
=
\sum_{a\in\mathcal A_p}\log F_a(p).
\label{eq:census-response}
\end{equation}

This census retains the individual active-modulus structure of each
exponent and is not assumed to depend on
\(p\) only through \(S(p)\). The resulting \(H(p)\) is the uncentered
\emph{order-selected structural score}. The identity product constraint
selects the family; it does not require every coordinate to equal one.
Its proposed transmission to primality is defined in
Section~\ref{sec:reference-anchor}.

The boundary \(a<p\) is structurally tied to the sector scale rather than
chosen to optimize the historical allocation.
Equations~\eqref{eq:psing-large-a} and~\eqref{eq:survival-an} imply,
uniformly for \(a\gg n\),
\begin{equation}
-\log s_{a,n}
=
\frac{n-1+\mathbf 1_{\{n=2\}}}{a-1}
+O\!\left(\frac{n^2}{a^2}\right).
\label{eq:large-a-filter}
\end{equation}
Averaging the factorwise relative response over~\eqref{eq:wn} as
in~\eqref{eq:survival-average}, and only then taking the negative
logarithm, gives

\begin{equation}
\log F_a(p)
=
\frac{\mu_p-1+w_2(p)}{a-1}
+O\!\left(\frac{\mathbb E_w[n^2]}{a^2}\right),
\qquad
\mu_p=\sum_n n\,w_n(p).
\label{eq:large-a-average}
\end{equation}

Since every supported sector satisfies
\(n\le N_{\max}(p)\asymp p/\log p\), the ratio \(n/a\) tends uniformly
to zero over \(a>p\) as \(p\) grows. Such moduli therefore lie
asymptotically in the weak-response regime of~\eqref{eq:large-a-average}. The exceptional case \(a=p\), which may
divide \(A_p\) for base-2 Wieferich exponents, is excluded
from the coprime CRT-channel construction.
This observation motivates the finite-scale cutoff \(a<p\) adopted
for the principal census, but does not establish
that the aggregate contribution of larger moduli
is negligible. The census and its subsequent
average-factor reduction are defined on the
active set \(a<p\); no approximation to the
complete contribution of all prime divisors
of \(A_p\) is asserted. Section~\ref{sec:uncentered-census} examines how
the resulting redistribution accumulates within
this selected range.

\section{Order-family census, population compression, and redistribution}
\label{sec:census-results}

\subsection{Detailed order-family census and target averaging}
\label{sec:uncentered-census}

For each of the \(5{,}761{,}450\) prime exponents \(13\le p<10^8\), the exact
active set~\eqref{eq:active-set} was enumerated through multiplicative-order
progressions and every local factor~\eqref{eq:Fa-explicit} was evaluated. The score \(H(p)\) is computed from these active moduli and the sector weights.
Neither \(\tau(p-1)\) nor \(S(p)\) serves as a summary predictor, and no
Mersenne-prime label enters the calculation.

As an unanchored diagnostic, one may insert the uncentered census response
directly into the same logarithmic mapping, without average-factor compression
or reference anchoring, and define
\begin{equation}
u_{\mathrm{cens}}(p)=\frac{\log p+H(p)}p.
\label{eq:census-kernel}
\end{equation}

The local bounds give the unanchored score an endpoint interpretation:
\[
\begin{gathered}
0\le H(p)\le\sum_{a\in\mathcal A_p}\log(a-1),\\
\log p\le\log p+H(p)
\le\log\!\left[p\prod_{a\in\mathcal A_p}(a-1)\right].
\end{gathered}
\]
Thus \(H(p)\) measures the effective auxiliary restriction assigned after
residue and sector averaging, on the additive log-density scale between
the neutral and fixed-channel bounds of the representation. These are
benchmark endpoints; their primality interpretation remains the separate
redistribution hypothesis of Section~\ref{sec:reference-anchor}.
For each interval, the higher-\(S\) allocation is its conditioned
case total multiplied by the normalized share of
\(u_{\mathrm{cens}}\) in its higher-\(S\) region; the six
allocations are then summed. Using the Wagstaff-balanced regions defined in
Appendix~\ref{app:wagstaff_balanced} and conditioning on the six observed
interval totals, this diagnostic gives
\[
28.610737/18.389263.
\]
This unanchored diagnostic measures score weighting across the fixed
\(S\) regions. It is distinct from the reference-anchored redistribution
of Section~\ref{sec:reference-anchor}. The idealized
\(23.5/23.5\) serves as the presentation baseline. The actual finite-partition Wagstaff mass share is
\(23.5544/23.4456\) (Appendix~\ref{app:order-control}); the exact
weighted-subset null mean is a different quantity,
\(\mathbb E_0X_+=23.327976\) (Appendix~\ref{app:wagstaff_balanced}).

\paragraph{Selection versus target-specific response.}
Keep the same exact active set \(\mathcal A_p\), and average the retained
response over targets before inversion:
\begin{equation}
H_{\mathrm{ord}}(p):=-\sum_{a\in\mathcal A_p}
\log\!\left[\sum_nw_n(p)s_{a,n}^{\mathrm{gen}}\right],
\qquad
H=H_{\mathrm{ord}}+\Delta_{\mathrm{tgt}}.
\label{eq:order-target-score}
\end{equation}
This control changes the target-dependent response, not the order-selected
family, sector measure, or aggregation. Its full-population allocation is
\(28.442566/18.557434\), against \(28.610737/18.389263\) for \(H\).
The comparison is label-independent apart from the conditioned interval totals.

For an interval \(I\) and its fixed higher region \(I_+\), put
\(B_I=\sum_I\log p/p\) and
\(\pi_{0,I}=\sum_{I_+}\log p/p\,/B_I\). For any score \(G\), define
\begin{equation}
D_I(G)=\sum_{p\in I_+}\frac{G(p)}p
-\pi_{0,I}\sum_{p\in I}\frac{G(p)}p,
\qquad
q_I(G)-\pi_{0,I}=
\frac{D_I(G)}{B_I+\sum_{p\in I}G(p)/p},
\label{eq:completion-contrast}
\end{equation}
Here \(q_I(G)\) is the higher-region share of \((\log p+G(p))/p\),
defined when these weights are positive on \(I\). The linear contrast
\(D_I(G)\) is also defined for signed score components.

Linearity gives \(D_I(H)=D_I(H_{\mathrm{ord}})+D_I(\Delta_{\mathrm{tgt}})\).
For additive attribution of the conditioned displacement, use the original
score's denominator for every component, where \(N_I\) is the observed
case count in interval \(I\):
\[
C(G):=\sum_I N_I\frac{D_I(G)}{B_I+\sum_{p\in I}H(p)/p}.
\]
Then \(C(H)=5.056296\), \(C(H_{\mathrm{ord}})=4.876617\), and
\(C(\Delta_{\mathrm{tgt}})=0.179679\). Thus target averaging preserves
\(96.45\%\) of this displacement from the finite-partition baseline
\(23.554441\). The percentage uses additive contrasts with common
denominators, not ratios of separately normalized allocations.

The target correction is substantial in \([13,10^3)\), contributing
\(0.096032\) of the contrast \(0.189318\). The aggregate percentage therefore
does not imply uniform target independence. Interval-level results,
modulus-count decomposition, and numerical bounds are given in Supplement~S.1.

The redistribution within the selected range is
concentrated at small moduli. The underlying census shows
that the fixed cutoff \(a<1000\) contains 31.43\% of all
active pairs and accounts for 92.82\% of \(\sum_p H(p)\).
Table~\ref{tab:census-cutoffs} shows that this cutoff captures 97.27\% of the
full unanchored higher-\(S\) displacement obtained with the \(a<p\)
census. In particular, extending the cutoff from \(a<1000\)
to \(a<p\) changes the conditioned allocation only from
\(28.4710/18.5290\) to \(28.6107/18.3893\).
The diagnostic stabilizes within this census as the cutoff increases;
this supplies neither a bound on the omitted \(a\geq p\) contribution nor
a causal ranking of moduli.

The \(a<N_{\max}(p)\) row is a separate moving-boundary diagnostic rather than
a further cumulative cutoff: depending on \(p\), it may exclude pairs counted
under \(a<1000\) and include pairs above that threshold, so no monotonicity
with the fixed-cutoff rows is expected.

\begin{table}[H]
\centering
\small
\begin{tabular}{lccc}
\toprule
Modulus selection & Higher \(S\) & Lower \(S\) & Share of full unanchored displacement \\
\midrule
Wagstaff 50/50 reference & 23.5000 & 23.5000 & 0\% \\
\(a<10\) & 24.5406 & 22.4594 & 20.36\% \\
\(a<30\) & 25.9405 & 21.0595 & 47.75\% \\
\(a<100\) & 27.2405 & 19.7595 & 73.19\% \\
\(a<300\) & 28.1076 & 18.8924 & 90.16\% \\
\(a<1000\) & 28.4710 & 18.5290 & 97.27\% \\
\(a<N_{\max}(p)\) & 28.2396 & 18.7604 & 92.74\% \\
\(a<p\) & 28.6107 & 18.3893 & 100\% \\
\bottomrule
\end{tabular}
\caption{Unanchored census diagnostics under cumulative modulus cutoffs. The
\(a<N_{\max}(p)\) boundary varies with \(p\) and is not nested uniformly with
the fixed numerical cutoffs. Percentages in the final column report the
fraction of the full \(a<p\) unanchored higher-\(S\) displacement, measured
from the idealized \(23.5/23.5\) presentation baseline.}
\label{tab:census-cutoffs}
\end{table}

\subsection{Population compression of the order-family score}

The count approximation~\eqref{eq:L-regression-new} describes how many active moduli occur, but not
their heterogeneous intensities. From the effective modular
aggregation~\eqref{eq:census-response}, we have the identity
\[
\log\lambda_{\mathrm{cens}}(p)
=\sum_{a\in\mathcal A_p}\log F_a(p)
=L(p)\log\bar F_{\mathrm{geom}}(p),
\]
where \(\bar F_{\mathrm{geom}}(p)\) is the geometric mean of the
active factors for that exponent.

Motivated by \(L(p)\simeq \kappa_I\tau(p-1)\), we compress this aggregate
logarithmic response onto the same structural scale. The
ordinary least-squares (OLS) coefficient through the origin over interval \(I\) is
\begin{equation}
\theta_I
=
\frac{\displaystyle\sum_{p\in I}\tau(p-1)H(p)}
{\displaystyle\sum_{p\in I}\tau(p-1)^2},
\qquad
H(p)\approx\theta_I\tau(p-1).
\label{eq:theta-projection}
\end{equation}

Defining the projected effective factor by
\(\bar F_I^{\mathrm{proj}}:=\exp(\theta_I/\kappa_I)\), the reduced
representation becomes
\begin{equation}
\begin{aligned}
\lambda_{\mathrm{red}}(p,S)
&=
\left(\bar F_I^{\mathrm{proj}}\right)^{\kappa_I(\log p)^S},
\\
H_{\mathrm{red}}(p,S)
&=\log\lambda_{\mathrm{red}}(p,S)
=\theta_I(\log p)^S.
\end{aligned}
\label{eq:FpowerL}
\end{equation}

The projected factor \(\bar F_I^{\mathrm{proj}}\) summarizes the
heterogeneous local factors already computed in the label-independent
census. It is obtained from \(H(p)\), not from \(L(p)\) alone, and is not
fitted to the historical Mersenne-prime split. The projection reintroduces
\(\tau(p-1)\) and supplies a label-independent parametrization of the
selected score for the redistribution hypothesis.
Section~\ref{sec:uncentered-census} gives its uncompressed diagnostic;
Section~\ref{sec:limitations} discusses the shared arithmetic of the score
and covariate.

For comparison, the direct geometric mean across all active
pairs in \(I\) is
\[
\log\bar F_I^{\mathrm{geom}}
=
\frac{\sum_{p\in I}H(p)}{\sum_{p\in I}L(p)}.
\]
The projected factor is used because its logarithmic coefficient
minimizes the squared error in approximating \(H(p)\) by
\(\theta_I\tau(p-1)\). The interval coefficients and the centered
coefficient of determination for the through-origin fit are shown in
Table~\ref{tab:effective-F}. For each interval \(I\), the latter is
\[
R_{H,I}^{2}
:=
1-
\frac{\sum_{p\in I}[H(p)-\theta_I\tau(p-1)]^2}
     {\sum_{p\in I}[H(p)-\overline H_I]^2},
\qquad
\overline H_I
:=
\frac{1}{|\mathcal P_I|}\sum_{p\in\mathcal P_I}H(p).
\]
Substituting this approximation into the unanchored diagnostic
\eqref{eq:census-kernel} gives
\[
28.585964/18.414036,
\]
only \(0.0248\) lower in the higher-\(S\) mass allocation than the complete
census. This comparison measures how closely the projection preserves the
conditioned allocation; it does not establish accuracy for other census
observables.

\begin{table}[H]
\centering
\small
\begin{tabular}{lrrr}
\toprule
Exponent interval & \(\bar F_I^{\mathrm{proj}}\) & \(\theta_I\) & \(R^2_{H,I}\) \\
\midrule
\([13,10^3)\) & 1.361163 & .161363 & .8103 \\
\([10^3,10^4)\) & 1.497133 & .234984 & .8246 \\
\([10^4,10^5)\) & 1.642612 & .308550 & .7999 \\
\([10^5,10^6)\) & 1.763574 & .390900 & .7923 \\
\([10^6,10^7)\) & 1.978867 & .490918 & .7950 \\
\([10^7,10^8)\) & 2.235156 & .602696 & .8020 \\
\midrule
Global & 2.213217 & .593255 & .7986 \\
\bottomrule
\end{tabular}
\caption{Population compression regressions
\(H\approx\theta_I\tau=\kappa_I\tau\log\bar F_I^{\mathrm{proj}}\), using the count
slopes in Table~\ref{tab:count-scaling}.}
\label{tab:effective-F}
\end{table}

The variation of \(\bar F_I^{\mathrm{proj}}\), from \(1.361\)
to \(2.235\), shows that a constant projected per-modulus
factor is not a faithful description over all six decades.
For the complete population, the corresponding label-independent
single-coefficient approximation is
\begin{equation}
\kappa=0.746752396760,
\qquad
\bar F^{\mathrm{proj}}=2.213216539032,
\qquad
\theta=\kappa\log\bar F^{\mathrm{proj}}=0.593255130346.
\label{eq:global-coefficients}
\end{equation}

\subsection{Identity assignments and boundary compatibility}
\label{sec:identity-boundary}

To compare the identity family with coherent shifted families, let
\(k\ge1\) be an integer and \(p>k\) prime. Define
the coherent shifted family and its realized target by
\[
\mathcal A_{p,k}
:=\{a<p:\ a\text{ odd prime},\ \operatorname{ord}_a(2)\mid p-k\},
\qquad t_k(a):=2^k-1\pmod a.
\]
The order relation gives \(M_p\equiv t_k(a)\pmod a\); for an active
modulus \(a<p\), this target is nonzero. Let \(\Omega_{a,n}(t)\) be
the set of nonzero residue tuples whose coordinate product is \(t\),
extending~\eqref{eq:omega-an}; put \(\mathbf1_n=(1,\ldots,1)\).

\begin{lemma}[Universal identity assignment]
\label{lem:universal-identity-assignment}
For \(a\in\mathcal A_{p,k}\) and \(n\ge1\),
\[
\mathbf1_n\in\Omega_{a,n}(t_k(a))
\quad\Longleftrightarrow\quad t_k(a)=1
\quad\Longleftrightarrow\quad a\mid2^{k-1}-1.
\]
Thus the assignment is available for every active modulus when \(k=1\).
For each fixed \(k>1\), infinitely many prime exponents have an active
modulus whose target fiber does not contain \(\mathbf1_n\).
\end{lemma}
\begin{proof}
The coordinate product is one, and \(a\) is odd, giving the equivalences.
Fix \(k>1\), choose a prime \(\ell>k\), and let \(a\) be a prime divisor
of \(2^\ell-1\). Its order divides the prime \(\ell\) and cannot be one,
so \(\operatorname{ord}_a(2)=\ell\). Since \((k,\ell)=1\), Dirichlet's
theorem gives infinitely many primes \(p\equiv k\pmod\ell\) with
\(p>\max\{a,k\}\).
Then \(a\in\mathcal A_{p,k}\), but identity availability would require
\(\ell\mid k-1\), impossible because \(\ell>k\).
\end{proof}

\begin{corollary}[Unique homogeneous representative across formal sectors]
Let \(G\) be a group and \(r,t\in G\). The same homogeneous assignment
\((r,\ldots,r)\) has fixed product target \(t\) for every integer
\(n\ge2\) if and only if \(r=t=1\). It is enough to require
\(r^n=r^{n+1}=t\) for two consecutive lengths.
\end{corollary}
\begin{proof}
If \(r^n=r^{n+1}=t\), multiplication by \(r^{-n}\) gives \(r=1\),
and hence \(t=1\). Conversely, \(r=t=1\) works for every length.
\end{proof}

The distinction underlying the selected family is primarily
between auxiliary moduli, not between the identity and
nonidentity residues within a selected modulus. For any odd
prime $a<p$, the classical bound $q\ge2p+1$ for prime divisors
of $M_p$ ensures that $M_p\bmod a\ne0$. Write
\[
 t_a(p):=M_p\bmod a\in\mathbb F_a^\times,
 \qquad
 \Omega_{a,n}(t_a)
 :=\left\{(u_1,\ldots,u_n)\in(\mathbb F_a^\times)^n:
          \prod_{i=1}^n u_i=t_a\right\}.
\]
The product target is fixed by $M_p$, independently of how
its factors are represented. The all-identity assignment
belongs to this fiber if and only if its target is one:
\[
 (1,\ldots,1)\in\Omega_{a,n}(t_a(p))
 \quad\Longleftrightarrow\quad
 t_a(p)=1
 \quad\Longleftrightarrow\quad
 a\in\mathcal A_p.
\]
This equivalence holds for every formal factorization length
$n\ge2$. It characterizes exactly the auxiliary moduli whose
product fibers admit the identity as a homogeneous
representative, independently of the formal factorization length.

Within a selected modulus, the individual residues may be
nonidentity residues, provided that their product is one.
The corollary shows, moreover, that the identity is the
unique \emph{fixed homogeneous} representative compatible
with product target one across all formal lengths. Thus the
selection has two levels: the product target selects the
eligible moduli, and homogeneous universality selects the
identity as their canonical representative. This uniqueness
does not apply to nonhomogeneous residue tuples or to
representatives chosen separately for each sector.

The set of all finite ordered residue sequences over
\(G_a=(\mathbb Z/a\mathbb Z)^\times\), including the empty sequence,
forms a monoid under concatenation. Among its fixed-product
families, only the product-one family is a submonoid:
concatenating two sequences with product \(t\) gives product
\(t^2\), which equals \(t\) only when \(t=1\). This closure applies
to arbitrary residue assignments, not only homogeneous ones,
and implies no additional primality restriction.

\paragraph{Localization across factor-number sectors.}
The same sector measure quantifies this homogeneous distinction.
With $w_*(p):=\max_{2\le n<R_p}w_n(p)$, define
\[
 W_a^{\mathrm{hom}}(p):=
 \max_{r\in\mathbb F_a^\times}
 \sum_{\substack{2\le n<R_p\\r^n=t_a(p)}}w_n(p).
\]
The identity target gives $W_a^{\mathrm{hom}}(p)=1$.
For a complementary target, $d_a:=\operatorname{ord}_a(t_a(p))\ge4$,
and unimodality of the sector weights gives
\[
 W_a^{\mathrm{hom}}(p)
 \le \frac1{d_a}+w_*(p)
 \le \frac14+w_*(p).
\]
Supplement~S.1 proves the bound and reports a finite-scale audit.
This diagnostic fixes one homogeneous residue across sectors;
it is distinct from the singleton response averaged over each full
fiber. It neither bounds the total complementary response nor
establishes that a residual primality multiplier is close to one.

The identity representative also has a conditional
factor-support interpretation. Together with the classical
condition $q_i\equiv1\pmod{2p}$, its genuine auxiliary
congruence $q_i\equiv1\pmod a$ implies
\[
 q_i\equiv1\pmod{2pa},\qquad q_i\ge 2pa+1.
\]
This statement concerns the all-identity configuration;
the product-one fiber does not force actual prime factors
to occupy that configuration.

Every fixed auxiliary residue channel has the same limiting
relative density $1/(a-1)$, which supplies the reference endpoint
used in the declared response. The effective response
$\rho_a=\overline s_a(p)$, however, is obtained by averaging the
declared response over the entire product-one fiber
(Sections~\ref{sec:local-survival} and
\ref{sec:effective-aggregation}), rather than imposing the
identity congruence on individual factors. Its inverse,
\[
 F_a(p)=\rho_a^{-1},\qquad
 -\log\rho_a=\log F_a(p),
\]
is a scalar response index, not an additional
arithmetic modulus.

For a modulus $a\notin\mathcal A_p$, the product target
$t_a(p)$ is not one. Consequently, its fiber has no
all-identity representative and is ineligible for the
\emph{same identity-channel representation}. Such a modulus may impose
other arithmetic restrictions or contribute to primality
through a different representation. The criterion neither
makes complementary moduli neutral nor excludes other
heuristic primality constructions. The separate question
of transmitting the selected effective response to
Mersenne primality is addressed in
Section~\ref{sec:reference-anchor}.

\paragraph{Boundary compatibility at \texorpdfstring{$a=p$}{a=p}.}
Although the effective census is restricted to $a<p$, the
formal order condition at the boundary $a=p$ provides an
additional distinction between the coherent shifted families.
For any odd prime $p>k$, Fermat's theorem gives
$2^{p-k}2^{k-1}\equiv1\pmod p$, and hence
\[
p\mid 2^{p-k}-1
\quad\Longleftrightarrow\quad
p\mid 2^{k-1}-1.
\]
Thus the $k=1$ family is boundary-compatible for every odd
prime exponent, whereas each fixed $k>1$ has at most
finitely many exceptional boundary-compatible primes.

Moreover, every prime factor $q$ of $M_p$ already satisfies
$q\equiv1\pmod{2p}$, so the identity residue at $a=p$
is compulsory in every realized factorization. The
universal boundary compatibility of the $k=1$ family
therefore links its product-one fibers to the classical
admissible-factor support. It does not, however, supply
an additional primality exclusion: the condition modulo
$p$ is already incorporated in Wagstaff's framework,
and $a=p$ is excluded from the coprime-CRT census.
No additional modular response factor is inferred
from this boundary observation.

\subsection{Reference anchoring and conditional redistribution}
\label{sec:reference-anchor}

Section~\ref{sec:structural-density-bridge} motivates a possible logarithmic
extension of Wagstaff's scale through conditional
factor support. Sections~\ref{sec:local-survival}
and~\ref{sec:effective-aggregation} construct
the corresponding declared response, while
Section~\ref{sec:identity-boundary} establishes
the algebraic role of the identity representative.
None of these results derives a primality law.

Every realized factorization of \(M_p\) satisfies its product congruences
modulo all odd auxiliary primes \(a<p\). These are consequences of the
factorization identity, not independent primality exclusions.
For an odd prime \(a<p\), the selected family is characterized exactly by
\[
a\in\mathcal A_p
\quad\Longleftrightarrow\quad M_p\equiv1\pmod a
\quad\Longleftrightarrow\quad \operatorname{ord}_a(2)\mid p-1.
\]
The selected family is exactly the family whose product
target admits the identity tuple in every formal sector.
Section~\ref{sec:identity-boundary} proves that no other
fixed homogeneous representative has this universal
product-one property. The family is characterized algebraically
by identity-channel compatibility, rather than defined through
a numerical ranking of competing shifted families.
Its relation to the divisors of $p-1$ follows from the
multiplicative-order criterion. Neither the admissibility
condition nor the identity's conditional factor-support
bound determines the declared response or its amplitude.

\paragraph{Log-density-equivalent scale.}
We represent a positive relative response on the
additive logarithmic scale
\[
J(\rho):=-\log\rho.
\]
For the classical factor support
$q\equiv1\pmod{2p}$, the reference density and
its logarithmic index are
\[
\rho_{\mathrm W}
=\frac1{\varphi(2p)}
=\frac1{p-1},
\qquad
J_{\mathrm W}
=\log(p-1)\sim\log p.
\]

For $B>1$ and $\rho>0$, define the formal
Mertens-equivalent scale
\[
B_{\mathrm{eq}}(B,\rho):=\frac{B}{\rho},
\qquad
\log B_{\mathrm{eq}}
=\log B-\log\rho.
\]
By Mertens' theorem, the product of $(1-1/\ell)^{-1}$ over
primes $\ell<B_{\mathrm{eq}}$ has leading term
$e^\gamma\log B_{\mathrm{eq}}$ as $B_{\mathrm{eq}}\to\infty$.
This equivalence parametrizes logarithmic
intensity only; it identifies neither a
finite sieve product nor an actual factor cutoff.

For the selected modular response,
Equation~\eqref{eq:census-response} gives
\[
-\log\rho_{\mathrm{cyc,cens}}(p)
=
\sum_{a\in\mathcal A_p}\log F_a(p)
=
H(p).
\]
At leading Wagstaff order, its unanchored
density-equivalent representation is therefore

\[
J_{\mathrm{sel}}(p):=\log\!\bigl[p\lambda_{\mathrm{cens}}(p)\bigr]
=\log p+\sum_{a\in\mathcal A_p}\log F_a(p)
=\log p+H(p).
\]
This is a definitional identity on the chosen equivalent scale, not a
prime-density theorem or a physical factor cutoff. The unanchored
continuation \(J_{\mathrm{sel}}(p)/p\) is the diagnostic in
Equation~\eqref{eq:census-kernel}; it is distinct from the principal
compressed and anchored rule below.

\paragraph{Restricted effective-channel continuation hypothesis.}
Wagstaff's marginal reference scale allows, in principle,
finite-scale conditional heterogeneity. We hypothesize that
the relative variation of the selected index, after
label-independent population compression and reference
anchoring, can represent such heterogeneity by
redistributing Wagstaff weight.

The restriction to \(\mathcal A_p\) is a selection principle,
not an auxiliary sieve over every modulus \(a<p\).
Precisely these moduli satisfy \(M_p\equiv1\pmod a\),
a condition indexed by divisors of \(p-1\), and their
product fibers admit the same homogeneous identity
representative for every formal factor-number sector.
The localization bound of
Section~\ref{sec:identity-boundary} and Supplement~S.1
quantifies the weaker cross-sector persistence of
complementary homogeneous channels. This difference in cross-sector persistence
motivates privileging the selected family, but does not
determine the amplitude of the singleton response,
which averages full product fibers.

The auxiliary product congruences for all \(a<p\) follow
from the same factorization identity and are not
additional independent primality exclusions.
Applying the declared response indiscriminately to
all such moduli would therefore define a different,
non-derived construction, not a complete sieve.
Neither this observation nor the homogeneous
localization result proves that unselected arithmetic
is neutral or already absorbed exactly in Wagstaff's
reference; it may reinforce, offset, or invalidate
the selected-index prediction. The compression,
reference anchoring and transmission to primality
remain explicit modeling choices, and adequacy as a
complete conditional primality model is unproved.

The population compression~\eqref{eq:global-coefficients} gives the
uncentered inverse-response index
\[
\lambda_{\mathrm{red}}(p,S)=\exp\!\left(\theta(\log p)^S\right).
\]
Define the relative factor at the reference \(S=1\) by
\begin{equation}
\Lambda(p,S):=
\frac{\lambda_{\mathrm{red}}(p,S)}{\lambda_{\mathrm{red}}(p,1)},
\qquad \Lambda(p,1)=1.
\label{eq:relative-anchor}
\end{equation}
The reference \(S=1\), equivalently \(\tau(p-1)=\log p\), is a
modeling choice rather than a population mean or median. It subtracts
the selected score at this declared reference before continuation on the
Wagstaff scale; neither CRT nor Wagstaff's derivation determines this anchor.

The \(S=1\) subtraction centers only the selected
structural score; it does not identify the population
mean of unselected arithmetic or prove that its
conditional contribution is constant. Any residual
modulation beyond the restricted continuation
remains unspecified.

Section~\ref{sec:historical-comparison} examines the alternative
reference \(\tau_0=A_0\log p\). The reference-centered score is
\begin{equation}
D(p,S):=\log\Lambda(p,S)
=\theta\left[(\log p)^S-\log p\right].
\label{eq:relative-log-response}
\end{equation}
The restricted continuation hypothesis now takes
the additive form
\begin{equation}
J_{\mathrm{eff}}(p,S)
=
J_{\mathrm W}(p)+\log\Lambda(p,S).
\label{eq:reduced-conditioned-kernel}
\end{equation}
This is the proposed transmission of the
reference-anchored equivalent response onto
Wagstaff's logarithmic scale, not a consequence
of the CRT construction or Dirichlet's theorem.

The subsequent factor \(\beta_I\) serves the distinct purpose of
preserving total Wagstaff mass on the finite calibration population.

At leading Wagstaff order, replacing \(J_{\mathrm W}=\log(p-1)\) by
\(\log p\) gives the proposed logarithmic rule
\begin{equation}
\log[p\Lambda(p,S)]
=\log p\left[(1-\theta)+\theta(\log p)^{S-1}\right].
\label{eq:convex-kernel}
\end{equation}
In this equivalent-scale parametrization, an unanchored score
\(H(p)\) maps a chosen classical reference scale \(B\) to \(Be^{H(p)}\).
After compression and \(S=1\) anchoring, the corresponding scale is
\((p-1)\Lambda(p,S)\), or \(p\Lambda(p,S)\) at leading Wagstaff order.
Because \(\Lambda\) is relative, this scale may lie above or below its
reference value; it is not asserted to be a lower bound on prime factors.

The proposed unnormalized exponent weight is
\begin{equation}
\widetilde w(p,S):=\frac{\log[p\Lambda(p,S)]}{p}
=\frac{(1-\theta)\log p+\theta(\log p)^S}{p}.
\label{eq:conditional-weighting-rule}
\end{equation}
Before calibration, the multiplier of \(b\) is \(1+D/\log p\),
not \(\Lambda\) itself. On the observed covariate,
\((\log p)^{S(p)}=\tau(p-1)\).
For a finite calibration population \(\mathcal P_I\), set
\(W_*(p):=\beta_I\widetilde w(p,S(p))\) and require
\begin{equation}
\sum_{p\in\mathcal P_I}W_*(p)
=\sum_{p\in\mathcal P_I}b(p).
\label{eq:marginal-consistency}
\end{equation}
Thus
\begin{equation}
\beta_I=
\frac{\displaystyle\sum_{p\in\mathcal P_I}\log p/p}
{\displaystyle\sum_{p\in\mathcal P_I}\log[p\Lambda(p,S(p))]/p}.
\label{eq:beta-definition}
\end{equation}
Conservation holds exactly on the calibration
population; global calibration need not preserve interval mass, and no
asymptotic equivalence is proved.
The construction redistributes the reference mass rather than creating
additional aggregate mass. Embedding these weights in Wagstaff's absolute
heuristic would multiply them by the common constant \(C\), which cancels
from conditional allocations. Interpreting the weights as conditional
primality densities remains the hypothesis; unresolved arithmetic dependence
beyond this rule is not modeled.

For a coefficient \(\theta\) fixed within \(I\), put
\[
x(p)=\frac{\tau(p-1)}{\log p},\qquad
\bar x_{b,I}:=
\frac{\sum_{p\in\mathcal P_I}b(p)x(p)}
{\sum_{p\in\mathcal P_I}b(p)}.
\]
The normalized redistribution is
\begin{equation}
\boxed{\frac{W_*(p)}{b(p)}
=\frac{1-\theta+\theta x(p)}
{1-\theta+\theta\bar x_{b,I}}.}
\label{eq:normalized-redistribution}
\end{equation}
For \(\theta>0\), weight increases relative to Wagstaff when
\(x(p)>\bar x_{b,I}\). Thus \(S=1\) anchors the score, whereas
\(\beta_I\) centers the redistribution over the calibration population.
At the formal reference, \(\widetilde w(p,1)=b(p)\); normalization gives
\(\beta_I b(p)\), not necessarily \(b(p)\).
The multiplier \((1-\theta)+\theta x\) is nonnegative for all \(x>0\)
exactly when \(0\le\theta\le1\). The global \(\theta=0.593255130346\)
and all interval coefficients satisfy this condition. Positivity establishes
an admissible weighting rule, not primality relevance.

For the sensitivity reference \(\tau_0=A_0\log p\),
\[
\frac{W_{*,A_0}(p)}{b(p)}
=\frac{1-\theta A_0+\theta x(p)}
{1-\theta A_0+\theta\bar x_{b,I}},
\]
provided the rule is positive on the calibration population.
The census and its compression determine \(\theta\) independently of
the subsequent choice of \(A_0\). Exact arithmetic selects the family;
the response, sector measure, aggregation, compression, and anchor determine
the redistribution amplitude.

The factor \(\beta_I\) is fixed after the weighting rule and does not fit
the higher-/lower-\(S\) split. It cancels from mass shares and the
weighted-subset law conditioned on each interval count. When the population
is understood, write \(\beta\). Calibration on \(13\le p<10^8\) gives
\(\beta=0.646960\). With global \(\theta\) fixed, \(\beta_I\) varies
from approximately \(0.651\) to \(0.643\) across the six ranges.
This describes approximate scale stability of marginal normalization,
without testing residual within-interval \(S\)-dependence. Interval-specific
\(\theta_I\) instead gives \(0.873\) to \(0.639\), reflecting compression drift.

\subsection{Historical allocations and model sensitivities}
\label{sec:historical-comparison}

Table~\ref{tab:new-balanced-allocation} separates the global rule from the
interval-specific compression. The global column uses one coefficient; the interval column follows
the detailed census more closely at each finite scale. These allocations
summarize relative population mass, not pointwise primality probabilities.

The displayed \(S_{50}\) values are rounded boundary summaries of the
population-defined half-Wagstaff-mass partitions in each interval. The exact
partitions, including deterministic tie breaking, are specified in
Appendix~\ref{app:wagstaff_balanced}.

For the model columns in Tables~\ref{tab:new-balanced-allocation}--\ref{tab:anchor-sensitivity}, an ``allocation'' is
obtained by multiplying each interval's normalized model mass share by
its conditioned interval total. The Wagstaff column in
Table~\ref{tab:new-balanced-allocation} instead
reports the idealized \(50/50\) presentation baseline associated with
the Wagstaff-balanced construction. These displayed allocations should
not be confused with the exact means of the without-replacement
weighted-subset law used in Appendix~\ref{app:wagstaff_balanced}.

\begin{table}[H]
\centering
\scriptsize
\begin{tabular}{lrrrrr}
\toprule
Exponent interval & \(S_{50}\) & Wagstaff & Interval \(\theta_I\) & Global \(\theta\) & Observed \\
\midrule
\([13,10^3)\) & 1.4479 & 5.0/5.0 & 5.424/4.576 & 6.078/3.922 & 7/3 \\
\([10^3,10^4)\) & 1.2535 & 4.0/4.0 & 4.691/3.309 & 5.369/2.631 & 5/3 \\
\([10^4,10^5)\) & 1.1770 & 3.0/3.0 & 3.722/2.278 & 4.154/1.846 & 4/2 \\
\([10^5,10^6)\) & 1.1162 & 2.5/2.5 & 3.286/1.714 & 3.548/1.452 & 4/1 \\
\([10^6,10^7)\) & 1.0959 & 2.5/2.5 & 3.485/1.515 & 3.618/1.382 & 3/2 \\
\([10^7,10^8)\) & 1.1027 & 6.5/6.5 & 9.544/3.456 & 9.513/3.487 & 9/4 \\
\midrule
Total & & 23.5/23.5 & 30.154/16.846 & 32.281/14.719 & 32/15 \\
\bottomrule
\end{tabular}
\caption{Historical allocations under the redistribution model,
conditioned on the six interval totals. The regions and compression
coefficients use no higher-/lower-\(S\) Mersenne labels; the continuation
was formulated after the historical association was known.}
\label{tab:new-balanced-allocation}
\end{table}

\begingroup
\setlength{\intextsep}{4pt plus 1pt minus 1pt}

Table~\ref{tab:global-sensitivity} compares several label-independent one-coefficient
reductions. The population OLS value is the equal-weight through-origin
projection in~\eqref{eq:theta-projection}, whereas the Wagstaff-weighted OLS gives
greater weight to exponents with larger baseline weight \(\log p/p\). The
geometric per-modulus mean provides a third compression. The final row retains
the population OLS coefficient and the \(S=1\) reference anchor, uses the
weight \([\log(a(p)p)+\theta(\tau(p-1)-\log p)]/p\), and evaluates it
on the partitions rebuilt with the refined Wagstaff weight in
Appendix~\ref{app:wagstaff_balanced}.

\begin{table}[H]
\centering
\small
\begin{tabular}{lrr}
\toprule
Global specification & \(\theta\) & Higher/lower \(S\) \\
\midrule
Population OLS; \(\log p\) baseline & .593255 & 32.281/14.719 \\
Wagstaff-weighted OLS; \(\log p\) baseline & .459832 & 30.904/16.096 \\
Geometric mean per active modulus; \(\log p\) baseline & .697019 & 33.209/13.791 \\
Population OLS; \(\log(a(p)p)\) numerator & .593255 & 31.519/15.481 \\
\bottomrule
\end{tabular}
\caption{Sensitivity of the historical model allocation. The first three rows
vary the compression of \(H(p)\) to a single coefficient \(\theta\); the last
holds \(\theta\) and the reference anchor fixed, changes the Wagstaff
numerator, and uses the refined-weight partitions. The coefficients and
partitions are computed without Mersenne-prime labels; the allocations condition on
the observed interval totals.}
\label{tab:global-sensitivity}
\end{table}

The anchor coefficient \(A_0\) is chosen separately from the projection
coefficient \(\theta\). On the observed covariate \(S=S(p)\), write
\begin{equation}
\tau_0(p)=A_0\log p,
\qquad
\log\Lambda_{A_0}(p,S(p))
=\theta\left[\tau(p-1)-A_0\log p\right].
\label{eq:anchor-sensitivity}
\end{equation}
The corresponding numerator is
\[
\log[p\Lambda_{A_0}]
=\left(1-\theta A_0\right)\log p+\theta\tau(p-1).
\]
Table~\ref{tab:anchor-sensitivity} varies \(A_0\) over a fixed window
around one while holding the label-independent global coefficient \(\theta\) fixed
and recalculating only the marginal factor \(\beta\).

\begin{table}[H]
\centering
\small
\begin{tabular}{rrrr}
\toprule
\(A_0\) & \(1-\theta A_0\) & Global \(\beta\) & Higher/lower \(S\) \\
\midrule
.75 & .555059 & .590317 & 31.518/15.482 \\
.90 & .466070 & .623047 & 31.959/15.041 \\
1.00 & .406745 & .646960 & 32.281/14.719 \\
1.10 & .347419 & .672782 & 32.628/14.372 \\
1.25 & .258431 & .715627 & 33.205/13.795 \\
\bottomrule
\end{tabular}
\caption{Sensitivity of the redistribution model to the reference
\(\tau_0=A_0\log p\). The structural coefficient and the higher-/lower-\(S\)
regions are unchanged throughout.}
\label{tab:anchor-sensitivity}
\end{table}

For the generalized reference \(\tau_0=A_0\log p\), the reduced multiplier is
\((1-\theta A_0)+\theta(\log p)^{S-1}\). For \(A_0>0\), uniform
non-negativity for all \(S\) is equivalent to \(0\le\theta\le1/A_0\).
Across the values in Table~\ref{tab:anchor-sensitivity}, the most restrictive
case, \(A_0=1.25\), requires \(\theta\le0.8\). The global value
\(\theta=0.593255\) and the alternative coefficients \(0.459832,0.697019\) in
Table~\ref{tab:global-sensitivity} all satisfy this bound. This is a consistency
property of the tested reductions.

\noindent\textbf{Retrospective agreement.}
The principal population-OLS compression and \(S=1\) anchor give
\(32.28/14.72\), against \(32/15\) observed. The coefficient \(\theta\) comes
from the label-independent census; \(\beta\) restores marginal mass
without fitting the split. The continuation was formulated after the
historical association was known (Section~\ref{sec:limitations}).
The high-\(S\) displacement persists across the sensitivities, but its
amplitude varies; changing \(A_0\) alone gives \(31.518\)--\(33.205\)
higher-region cases. Neither the compression nor the anchor is uniquely
identified.

\endgroup

\section{Extrapolation beyond the historical range}
\label{sec:prospective}

The conditional weighting rule of Section~\ref{sec:reference-anchor} can be
evaluated beyond the historical calibration range without using
Mersenne-prime outcomes in the new interval. For \(10^8\le p<10^9\),
the fixed higher-\(S\) region is defined by \(S(p)\ge S_{50}\), with the
rounded Wagstaff-weighted median
\[
S_{50}=1.071704.
\]
The population contains \(45{,}086{,}079\) prime exponents, and its baseline
higher-region weight share is \(0.49999965\). Keeping the global coefficient
of Equation~\eqref{eq:global-coefficients} and the reference \(A_0=1\)
gives the model allocation
\begin{equation}
73.8841\%/26.1159\%.
\label{eq:prospective-share}
\end{equation}
Continuing the final historical interval coefficient \(\theta_6=0.602696\)
gives \(74.1270\%/25.8730\%\); this is a sensitivity calculation, not an
alternative to select after future observations. The result is the share
implied by the declared weighting rule on the complete exponent population.
The threshold and principal weighting specification are fixed as reported
in this dated manuscript version; subsequent comparisons should retain them
without outcome-dependent retuning.

For orientation only, the coarse Wagstaff decade expectation is
\[
\frac{e^\gamma}{\log2}\log10=5.9166,
\]
corresponding to illustrative allocations of approximately \(4.371/1.545\)
under the proposed redistribution. These absolute counts inherit the
uncertainty of the Wagstaff heuristic. The exponent \(136{,}279{,}841\) was
known before the extrapolation was specified and supplies no independent
confirmation. The discriminating power of a conditional comparison depends
on the number of cases and their allocation, as illustrated in Supplement~S.2.

The conditional comparisons in Supplement~S.2 quantify how distinct
weighting rules would distribute cases in a specified, adequately
ascertained eligible population. The \(73.88\%\) figure pertains to the
entire prime-exponent interval and does not automatically apply to a subset
selected for computational testing; its reference mass shares would have
to be recalculated for that subset. Comparisons based only on reported
discoveries, without reconstructing search eligibility and exposure, remain
descriptive. Future evidence could inform the persistence of the historical
association and the relative adequacy of the proposed redistribution;
Section~\ref{sec:limitations} discusses the limits of a mechanistic interpretation.

\section{Limitations and interpretation}
\label{sec:limitations}

The statistical association and proposed redistribution concern different
claims. Local analyses reuse historical cases and have overlapping windows,
ties, and boundary truncation; their inferential values are descriptive or
nominal. The population reference is exact for the specified weighted-subset
law and fixed partition. It conditions on only \(47\) cases and does not
derive the historical discovery process or account for all exploratory
choices that preceded the specification of \(S\).

The post hoc \(v_2(p-1)\)-conditioned design in
Appendix~\ref{app:order-control} finds no separately detectable odd-part
association. Under the same partition, the exact higher-region means are
\(23.42\) for Wagstaff, \(27.64\) for the unanchored census, and \(31.48\)
for the global reduction, against \(25\) observed. Write \(p-1=2^v u\)
with \(u\) odd. Since \(\tau(p-1)=(v+1)\tau(u)\), the component
projection introduces odd-part variation even into the contribution
from moduli whose order is a pure power of two. This diagnoses a loss of conditional fidelity,
but neither explains the whole discrepancy nor isolates compression from
anchoring. The order-class contributions are not a causal attribution
of the observed cases.

A separate post hoc threshold diagnostic finds $v_2(p-1)\ge3$ in
19 of the 47 historical cases, against an exact Wagstaff-conditional
expectation of 11.36 (upper-tail probability $0.009183$); this correlated
exploratory comparison is not an independent confirmation of the
full-divisor association.

An additional post hoc audit, reported in Supplement~S.3, tests
direct transmission of the uncompressed singleton score restricted to moduli with even multiplicative order dividing \(p-1\). After within-interval centering under Wagstaff weights, its retrospective higher-$S$ mass allocation is 30.11
of 47 cases, against 32 observed, without regressing $H$ onto
$\tau(p-1)$. Approximately 84\% of this modeled higher-$S$
displacement comes from mixed even orders, that is, orders with a nontrivial odd factor. After conditioning on $v_2(p-1)$, the raw and centered rules
give conditional full-set log-likelihood ratios of $+0.9608$ and
$-0.1623$, respectively, relative to Wagstaff. Neither comparison
provides convincing evidence of additional primality-related information
beyond the $2$-adic valuation.

This conditional limitation does not by itself contradict the
order-family interpretation. Indeed, activation of an even-order
modulus with \(\operatorname{ord}_a(2)=2^m t\), where \(t\) is odd,
requires \(m\leq v_2(p-1)\)
(Appendix~\ref{app:order-control}). The centered direct even-order
benchmark also gives a mass-share allocation of \(29.98\) cases to
the higher region of \(S_{2\mathrm{ad}}(p)\), obtained by replacing
\(\tau(p-1)\) with \(1+v_2(p-1)\) in \(S(p)\), against \(31\) observed
(Supplement~S.3),
without projection onto \(\tau(p-1)\). Thus, a predominantly
\(2\)-adically aligned redistribution is compatible with the
selected-order architecture. This compatibility neither establishes
transmission to primality nor distinguishes the modular mechanism
from a simpler model based only on \(v_2(p-1)\).

The cutoff \(a<p\) selects a subsystem; it does not bound the total omitted
contribution. The scale drift of \(\theta_I\), sector-prior sensitivity, and
anchor sensitivity limit amplitude identification. The historical global
allocation \(32.28/14.72\), against \(32/15\) observed, represents
specification-dependent retrospective agreement in direction and
magnitude rather than independent validation. Near-constant \(\beta_I\)
describes marginal normalization, not absence of residual within-interval
\(S\)-dependence.

The restricted continuation specifies selected-benchmark weights, without
determining the contribution of omitted arithmetic to a true conditional
primality law. If that law differs by a positive residual multiplier,
its normalized higher-\(S\) share depends on the covariance of this
multiplier with the higher-region indicator under the selected weights.
The index definition imposes no condition making that covariance small.

A separate post hoc sensitivity analysis examines
the effect of multiplying the restricted rule by
a positive residual proxy constructed from
shifted-family responses. The principal higher-\(S\)
allocation of 32.281 cases varies between 31.438
and 33.093 across the tested perturbations
\(\eta=-2,-1,0,1,2\), and decreases to 27.310
under the deliberately strong stress \(\eta=11\).
These values are illustrative sensitivity scenarios,
not confidence bounds or estimates of the
unresolved arithmetic contribution. They show
that the selected-index prediction may be altered
by residual structure not constrained by the
identity-channel criterion (Supplement~S.4).

The distinction between exact results and modeling choices is as follows.
\begin{itemize}
\item The activation criterion, CRT product reference, singleton formula,
and uniqueness of the homogeneous identity representative are exact.
The uniqueness statement concerns formal sector lengths. The density
ratio \(1/(a-1)\) holds in the fixed-modulus limit, without a uniform
finite-scale estimate over growing moduli.
\item The factorwise response in Section~\ref{sec:local-survival} is
declared. Equation~\eqref{eq:survival-an} is its exact expectation, not a
derivation of a factorization survival law. Its positive inverse is a
scalar equivalent index, not an arithmetic modulus or a prime-density law.
\item The sector weights, multiplicative response, and order of sector
marginalization and modular aggregation are benchmark choices.
Appendices~\ref{app:allocation-audit} and~\ref{app:sector-prior-sensitivity}
audit quantitative sensitivity to these choices.
\item The projection~\eqref{eq:FpowerL} is a population compression of the
selected score. Both \(H(p)\) and \(S(p)\) use the arithmetic of \(p-1\);
the projection's accuracy is therefore not independent evidence for the
primality relevance of \(S\).
\item The anchor~\eqref{eq:relative-anchor}, logarithmic continuation
\eqref{eq:reduced-conditioned-kernel}, and exponent weight
\eqref{eq:conditional-weighting-rule} specify the redistribution
hypothesis. Equation~\eqref{eq:beta-definition} preserves marginal
Wagstaff mass on the calibration population. This imposed normalization
cancels from the conditional comparison and does not validate transmission.
\end{itemize}

\subsection*{Family selection, local targets, and scope}
Target averaging within the same order-selected family preserves
96.45\% of the common-denominator benchmark displacement
(Supplement~S.1). This replaces the local target-specific response
while keeping the active set fixed. Thus the universal identity
representative motivates selection but does not by itself explain the
benchmark amplitude; active-modulus selection supplies most of this
constructed displacement. The all-moduli extension is a scope
control, not a derived complete sieve. Its near-neutral allocation
does not demonstrate cancellation or neutrality of complementary
arithmetic effects. Generic modulus-size and sector geometry largely
reproduce its negative contrasts.

The extrapolation in Section~\ref{sec:prospective} is a fixed numerical
implication of the chosen rule. Future observations could help assess
persistence, but even persistent enrichment would remain compatible with
a correlated arithmetic proxy rather than the proposed mechanism.

Another arithmetic variable correlated with \(S(p)\) could account for
the observed enrichment. The algebraic selection of the family does
not establish unique primality relevance.

\section{Conclusion}

Among the \(47\) known Mersenne-prime exponents \(13\le p<10^8\),
\(32\) lie in the higher-\(S\) regions of six Wagstaff-balanced
intervals. The historical allocation is unusual under the specified
conditional weighted-subset reference, with exact two-sided value
\(p_{\mathrm{two}}=0.01154\). The statistic and related analyses were
selected post hoc, so this nominal result does not provide confirmatory
evidence of persistent conditional heterogeneity. Conditioning on \(v_2(p-1)\) leaves no separately
detectable odd-part association.

The arithmetic framework identifies an order family linked exactly
to the cyclotomic structure of \(2^{p-1}-1\). The identity is the unique
homogeneous representative compatible with a fixed target at every formal
sector length. This motivates the selected family without deriving its
declared response. Residue and sector averaging define a formal
density-equivalent index, whose logarithmic aggregate is \(H(p)\).
Its population compression summarizes the benchmark reproducibly;
its alignment with \(S(p)\) also reflects their shared arithmetic.

Under the restricted effective-channel continuation hypothesis,
population compression and \(S=1\) anchoring yield a candidate conditional
refinement that preserves Wagstaff mass on the finite calibration
population. The retrospective allocation is \(32.28/14.72\), compared
with \(32/15\) observed. This agreement is specification-dependent and
was obtained after the historical association was known. The selected
index is defined without a neutrality assumption for complementary
arithmetic; its adequacy as a complete conditional primality model remains
unproved. The available sample has not demonstrated detectable
primality-related information beyond the \(2\)-adic valuation. A derivation
linking the modeled response to Mersenne primality remains missing;
the restricted transmission hypothesis is neither established nor ruled out.

The unchanged rule assigns \(73.88\%\) of its conditional weight
to the higher-\(S\) region in \(10^8\le p<10^9\). Its discriminating
power depends on the number and allocation of subsequent eligible
cases. Assessing whether the historical association persists and
deriving a principled link between the selected modular constraints
and Mersenne primality remain separate open problems.

\begin{appendices}

\section{Statistical analysis and robustness}
\label{appendix:statistical}

For each known Mersenne-prime exponent \(p\), we compare the divisor structure of \(p-1\) with that of nearby prime exponents of comparable size.

Let \(p_{(k)}\) denote the \(k\)-th prime.
For a given exponent \(p=p_{(i)}\) and a window parameter \(W\), define the control set

\[
\mathcal P_W(p)
=
\{p_{(i-j)}: 1\le j\le W_-(p)\}
\cup
\{p_{(i+j)}: 1\le j\le W\},
\]

where

\[
W_-(p)=\min(W,i-1).
\]

Thus the control set contains the next \(W\) primes after \(p\)
and up to \(W\) primes before \(p\), excluding \(p\) itself.

We define the percentile rank using mid-rank assignment:

\[
\pi_W(p)
=
\frac{
\#\{\,q\in\mathcal P_W(p): S(q)<S(p)\}
+\tfrac12\,\#\{\,q\in\mathcal P_W(p): S(q)=S(p)\}
}{
|\mathcal P_W(p)|
}.
\]

Thus the reported local percentiles rank the normalized
statistic \(S\), rather than the unnormalized divisor count
\(\tau(p-1)\). This definition uses mid-rank assignment to account for ties;
its inferential interpretation is qualified in
Appendix~\ref{app:local-summary}.

\subsection{Mersenne-prime exponent data}
\label{app:exponent-data}

Table~\ref{tab:Sp_percentiles} lists the Mersenne-prime exponents
\(p\ge13\) in the GIMPS dataset accessed September 2026, together with
their normalized divisor statistic \(S(p)\) and local percentile ranks.

Exponents below \(p = 13\) are excluded because the normalizing factor
\(\log\log p\) varies rapidly there.

\setlength{\LTpre}{6pt}
\setlength{\LTpost}{6pt}

\begingroup
\renewcommand{\arraystretch}{0.97}
\begin{longtable}{rrrrr}
\caption{Divisor parameter \(S(p)\) and local percentile ranks for the
Mersenne-prime exponents in the GIMPS dataset \cite{GIMPS}. Percentiles use
the next \(W\) primes after \(p\) and up to \(W\) before \(p\), excluding \(p\);
the left side is truncated near the lower boundary.}
\label{tab:Sp_percentiles}\\

\toprule
\(p\) & \(\tau(p-1)\) & \(S(p)\) & \(\pi_{1000}\) & \(\pi_{5000}\) \\
\midrule
\endfirsthead

\multicolumn{5}{l}{\small Table \ref{tab:Sp_percentiles} (continued).}\\
\toprule
\(p\) & \(\tau(p-1)\) & \(S(p)\) & \(\pi_{1000}\) & \(\pi_{5000}\) \\
\midrule
\endhead

\midrule
\multicolumn{5}{r}{\small Continued on next page.}\\
\endfoot

\bottomrule
\endlastfoot
13 & 6 & 1.9022 & 0.997 & 0.997 \\
17 & 5 & 1.5454 & 0.832 & 0.869 \\
19 & 6 & 1.6592 & 0.929 & 0.940 \\
31 & 8 & 1.6855 & 0.947 & 0.954 \\
61 & 12 & 1.7578 & 0.976 & 0.978 \\
89 & 8 & 1.3849 & 0.666 & 0.727 \\
107 & 4 & 0.8992 & 0.130 & 0.194 \\
127 & 12 & 1.5749 & 0.850 & 0.892 \\
521 & 16 & 1.5122 & 0.806 & 0.844 \\
607 & 8 & 1.1194 & 0.349 & 0.408 \\
1279 & 12 & 1.2629 & 0.512 & 0.599 \\
2203 & 8 & 1.0189 & 0.288 & 0.296 \\
2281 & 32 & 1.6944 & 0.950 & 0.958 \\
3217 & 20 & 1.4341 & 0.721 & 0.767 \\
4253 & 6 & 0.8440 & 0.114 & 0.107 \\
4423 & 16 & 1.3032 & 0.595 & 0.627 \\
9689 & 16 & 1.2507 & 0.558 & 0.587 \\
9941 & 24 & 1.4318 & 0.760 & 0.764 \\
11213 & 6 & 0.8025 & 0.096 & 0.100 \\
19937 & 24 & 1.3862 & 0.722 & 0.734 \\
21701 & 36 & 1.5573 & 0.887 & 0.882 \\
23209 & 16 & 1.2014 & 0.528 & 0.541 \\
44497 & 40 & 1.5561 & 0.892 & 0.882 \\
86243 & 16 & 1.1407 & 0.469 & 0.479 \\
110503 & 24 & 1.2960 & 0.664 & 0.669 \\
132049 & 60 & 1.6594 & 0.941 & 0.946 \\
216091 & 60 & 1.6323 & 0.935 & 0.936 \\
756839 & 8 & 0.7981 & 0.166 & 0.153 \\
859433 & 32 & 1.3254 & 0.735 & 0.727 \\
1257787 & 12 & 0.9404 & 0.288 & 0.280 \\
1398269 & 6 & 0.6762 & 0.060 & 0.061 \\
2976221 & 24 & 1.1763 & 0.567 & 0.574 \\
3021377 & 28 & 1.2329 & 0.633 & 0.642 \\
6972593 & 40 & 1.3379 & 0.767 & 0.770 \\
13466917 & 36 & 1.2806 & 0.748 & 0.736 \\
20996011 & 180 & 1.8383 & 0.991 & 0.990 \\
24036583 & 8 & 0.7340 & 0.129 & 0.128 \\
25964951 & 48 & 1.3643 & 0.801 & 0.798 \\
30402457 & 128 & 1.7045 & 0.974 & 0.973 \\
32582657 & 44 & 1.3275 & 0.745 & 0.753 \\
37156667 & 16 & 0.9700 & 0.354 & 0.351 \\
42643801 & 192 & 1.8344 & 0.993 & 0.992 \\
43112609 & 48 & 1.3504 & 0.803 & 0.789 \\
57885161 & 64 & 1.4424 & 0.865 & 0.868 \\
74207281 & 80 & 1.5125 & 0.918 & 0.909 \\
77232917 & 12 & 0.8571 & 0.220 & 0.227 \\
82589933 & 12 & 0.8560 & 0.221 & 0.221 \\
136279841 & 24 & 1.0846 & 0.498 & 0.501 \\

\end{longtable}
\endgroup

\subsection{Summary statistics and local reference comparisons}
\label{app:local-summary}

Summary statistics for the divisor statistic and the local percentile ranks
are presented in Table~\ref{tab:summary_stats}.

\begin{table}[ht]
\centering
\caption{Summary statistics for the divisor statistic and the local percentile ranks.
The reported \texttt{std} uses denominator \(n\), and the rows labelled \texttt{p10}, \texttt{p25}, \texttt{p75}, and \texttt{p90} use the nearest-observation convention.}
\label{tab:summary_stats}
\begin{tabular}{lrrr}
\toprule
Statistic & \(S(p)\) & \(\pi_{1000}\) & \(\pi_{5000}\) \\
\midrule
n & 48 & 48 & 48 \\
mean & 1.3158 & 0.6372 & 0.6484 \\
median & 1.3327 & 0.7285 & 0.7350 \\
std & 0.3174 & 0.2907 & 0.2891 \\
min & 0.6762 & 0.0600 & 0.0612 \\
p10 & 0.8560 & 0.1660 & 0.1944 \\
p25 & 1.1194 & 0.4690 & 0.4791 \\
p75 & 1.5561 & 0.8875 & 0.8824 \\
p90 & 1.6944 & 0.9499 & 0.9580 \\
max & 1.9022 & 0.9970 & 0.9972 \\
\bottomrule
\end{tabular}
\end{table}

Under an idealized uniform within-choice-set reference with distinct ranks,
the probabilities of falling above and below \(0.5\) are symmetric, apart
from a possible central rank. Ties, overlapping choice sets, and boundary
truncation limit the interpretation of this reference for the present data.
In particular, the earliest boundary-truncated sets span a wider exponent
range, over which the size-dependent Wagstaff baseline may vary appreciably.
The sign, Wilcoxon, and KS calculations below are therefore retained as
descriptive local summaries, rather than exact tests under the
Wagstaff-weighted null.

Using the sign statistic on the \(\pi_{5000}\) values, we observe
35 positive deviations out of 48 comparisons. The corresponding
binomial reference calculation,
\[
X\sim\mathrm{Bin}(48,1/2),
\]
gives
\[
p=0.0021
\qquad\text{(two-sided)}.
\]

The Wilcoxon signed-rank statistic applied to the centered values
\(\pi_{5000}-0.5\) gives
\[
W_+=893,
\]
with the nominal two-sided reference value
\[
p=0.0014.
\]

A descriptive Kolmogorov--Smirnov comparison with the
\(\mathrm{Uniform}(0,1)\) distribution gives
\[
D=0.290,
\qquad
p=0.00046.
\]

Because the percentile ranks are discrete and the local windows overlap,
the KS result is reported as a descriptive distributional comparison.
The stratified permutation procedure described below provides a
complementary conditional assessment.

\noindent
\textbf{Effect size (mean shift).}

We compare the mean \(S(p)\) of the Mersenne-prime exponents with the pooled
local-control distribution. Every control occurrence is retained, including
repeated occurrences in overlapping windows. The effect sizes describe
displacement relative to this pooled distribution.

Using the \(W=1000\) and \(W=5000\) control windows, we obtain

\[
\overline S_M=1.3158,
\qquad
\overline S_C^{(1000)}=1.1526,
\qquad
\overline S_C^{(5000)}=1.1383,
\]

yielding mean shifts

\[
\Delta^{(1000)}=0.1632,
\qquad
\Delta^{(5000)}=0.1774.
\]

The corresponding standardized mean differences use the root-mean-square
of the Mersenne and pooled-control sample standard deviations,
\(s_{\mathrm{av}}=[(s_M^2+s_C^2)/2]^{1/2}\). They are

\[
d^{(1000)}=0.510,
\qquad
d^{(5000)}=0.563.
\]

\subsection{Stratified conditional model}

To quantify the local association while preserving the local choice-set
structure, we fit a stratified conditional logistic regression
model~\cite{McFadden}.

Each choice set \(\mathcal C_i\) consists of one observed Mersenne-prime exponent
and its local prime controls:
\[
\mathcal C_i
=
\{p_i\}
\cup
\mathcal P_W(p_i),
\]
where \(\mathcal P_W(p_i)\) contains the \(W\) nearest primes below and the \(W\) nearest primes above \(p_i\)
(truncated near the lower boundary).

The conditional selection probability is modeled as
\[
\mathbb P(p \text{ selected in } \mathcal C_i)
=
\frac{\exp(\gamma_{\mathrm{loc}} S(p))}
{\sum_{q\in\mathcal C_i} \exp(\gamma_{\mathrm{loc}} S(q))}.
\]

For \(W=1000\), we obtain
\[
\widehat\gamma_{\mathrm{loc}}=0.66,\qquad p_{\rm LR}\approx0.012.
\]
For \(W=5000\), we obtain
\[
\widehat\gamma_{\mathrm{loc}}=0.84,\qquad p_{\rm LR}\approx0.003.
\]

At \(\gamma_{\mathrm{loc}}=0\), the model gives uniform within-set selection.
Window overlap and boundary truncation limit the likelihood-ratio values
to nominal local references; the size-dependent Wagstaff baseline and these
limitations are discussed in Appendices~\ref{app:local-robustness}
and~\ref{app:wagstaff_balanced}.

\subsection{Stratified permutation test}

To obtain a finite-sample descriptive reference under the adopted uniform
within-stratum scheme, we performed a stratified permutation test.

Within each choice set \(\mathcal C_i\), one prime was selected uniformly at random,
preserving the stratified structure of the conditional model.
For each permutation, we computed the centered stratified score statistic

\[
T
=
\sum_{i}
\left(
S(\tilde p_i)
-
\bar S_i
\right),
\qquad
\bar S_i
=
\frac{1}{|\mathcal C_i|}
\sum_{q\in\mathcal C_i} S(q),
\]

where \(\tilde p_i\) denotes the randomly selected prime in stratum \(i\).

The reference is the uniform within-stratum scheme, subject to the
limitations discussed in Appendix~\ref{app:local-robustness}.

Using \(20{,}000\) stratified permutations, we obtain the following
two-sided reference values:

\[
p_{\mathrm{perm}} = 0.0039 \quad (W=1000),
\qquad
p_{\mathrm{perm}} = 0.0012 \quad (W=5000).
\]
The reference implementation initializes the same fixed generator seed,
\(260300166\), separately for each window size.

Two-sided permutation \(p\)-values were computed as

\[
p
=
\frac{1+\#\{|T_{\mathrm{perm}}| \ge |T_{\mathrm{obs}}|\}}
{1+N_{\mathrm{perm}}},
\]

using the standard finite-sample correction for Monte Carlo permutation
tests.

These permutation values summarize the local enrichment under the adopted
comparison scheme. The complete-population Wagstaff-weighted randomization in
Appendix~\ref{app:wagstaff_balanced} provides the primary inferential
assessment under its specified conditional null.

\subsection{Distinct prime factors:
\texorpdfstring{\(\omega(p-1)\)}{omega(p-1)}}

To separate the contribution of the number of distinct prime factors
from divisor multiplicities, we consider the arithmetic function
\(\omega(n)\), which counts the number of distinct prime divisors of \(n\).

We define the normalized parameter

\[
S_\omega(p)=\frac{\log \omega(p-1)}{\log\log p}.
\]

The same percentile and permutation procedures give weaker enrichment
for \(S_\omega(p)\) than for \(S(p)\), with a small standardized effect.
This comparison alone does not isolate an independent multiplicity effect.

\subsection{Multiplicity-sensitive divisor ratio}
To construct a multiplicity-sensitive contrast in the prime
factorization of \(p-1\), we consider the ratio

\[
R(p)=\frac{\tau(p-1)}{2^{\omega(p-1)}}.
\]

The denominator \(2^{\omega(p-1)}\) corresponds to the divisor count of a
squarefree integer with the same number of distinct prime factors.
Thus \(R(p)\) measures the excess divisor structure arising from
prime-factor multiplicities.

We define the normalized parameter

\[
S_R(p)=\frac{\log R(p)}{\log\log p}.
\]

Repeating the statistical analysis with \(S_R(p)\) yields
enrichment comparable to that obtained for \(S(p)\), both in mean
percentile and in standardized effect size.

Table~\ref{tab:structural_comparison} summarizes the behavior of the
three structural parameters. While \(S(p)\) and \(S_R(p)\) show
moderate positive enrichment, \(S_\omega(p)\) exhibits only a weak
effect.

We retain \(S(p)\) as the primary structural statistic because the
divisors of \(p-1\) index the cyclotomic layers described in
Section~\ref{sec:algebraic}. The alternatives \(S_\omega(p)\) and \(S_R(p)\)
are post hoc robustness controls; the comparable enrichment in \(S_R(p)\)
does not motivate replacing the principal statistic.

\begin{table}[ht]
\centering
\begin{tabular}{lccc}
\toprule
Parameter & Structural component & Mean percentile & Effect size \(d\) \\
\midrule
\(S(p)\) & total divisor structure & 0.65 & 0.56 \\
\(S_\omega(p)\) & distinct prime factors & 0.58 & 0.28 \\
\(S_R(p)\) & prime-factor multiplicities & 0.65 & 0.53 \\
\bottomrule
\end{tabular}
\caption{Comparison of structural parameters based on the divisor
structure of \(p-1\), using the \(W=5000\) local controls. All effect
sizes use the same root-mean-square group-standard-deviation convention.}
\label{tab:structural_comparison}
\end{table}

\FloatBarrier

\subsection{Robustness analysis}
\label{app:local-robustness}

Overlapping windows induce dependence, and the asymmetric choice sets near
the lower boundary span wider exponent ranges. In those sets, uniform
within-stratum exchangeability need not approximate the size-dependent
Wagstaff baseline. The local \(p\)-values are descriptive; the complete
population analysis in Appendix~\ref{app:wagstaff_balanced} incorporates
Wagstaff weights and provides the primary assessment under its specified
conditional null.

For a fixed exponent and control window, the percentile is determined
arithmetically; this alone does not control selection of the historical
sample. The GIMPS milestones record reports first-test coverage below
141 million by July 24, 2026, covering all listed exponents
\cite{GIMPSMilestones}. With complete first-test coverage and correct
final classifications, prescreening changes how composites are eliminated,
rather than their classification. The small sample and exploratory
choice of summaries still limit inference.

The local procedures probe different summaries of the same historical sample
and should not be regarded as independent replications. Their agreement shows
only that the observed enrichment is not tied to a particular local statistic
or inferential procedure.

\subsection{Population-level Wagstaff-balanced randomization}
\label{app:wagstaff_balanced}

The local-window analyses above use overlapping comparison sets. We therefore
also use a population-level conditional design over all \(5{,}761{,}450\)
prime exponents \(13\le p<10^8\). The six disjoint exponent intervals are the
ones in Table~\ref{tab:new-balanced-allocation}, and we condition on their
observed Mersenne-prime counts
\[
(N_1,\ldots,N_6)=(10,8,6,5,5,13).
\]
This removes the distribution of discoveries across exponent scales from the
comparison and tests only the within-interval allocation by \(S\).

Let \(\mathcal P_b\) be the prime-exponent population in interval \(b\).
For each interval, sort its complete population by decreasing
\(S(p)\), breaking ties by decreasing \(p\), and take the prefix whose
cumulative baseline weight is closest to one
half of the interval total. This defines the higher-\(S\) region
\(\mathcal P_b^+\); its complement is \(\mathcal P_b^-\). The primary weight
is
\[
w_0(p)=\frac{\log p}{p}.
\]
The construction uses a prefix rather than assigning every boundary tie
to the same region. The thresholds reported in Table~\ref{tab:new-balanced-allocation}
are therefore boundary summaries, while the stored Boolean partition is the
exact object used in every calculation. No Mersenne-prime label enters the
partition.

Within interval \(b\), a subset \(A_b\subset\mathcal P_b\) of size \(N_b\) is
assigned the conditional product-weight probability
\begin{equation}
\Pr_0(A_b\mid |A_b|=N_b)
=
\frac{\displaystyle\prod_{p\in A_b}w_0(p)}
{\displaystyle
 \sum_{\substack{B\subset\mathcal P_b\\|B|=N_b}}
 \prod_{q\in B}w_0(q)}.
\label{eq:weighted-subset-null}
\end{equation}
This is Fisher's noncentral hypergeometric subset law with item-specific
weights. It incorporates the complete finite population, sampling without
replacement, the variation in \(\log p/p\), and the small residual imbalance
of the discrete half-weight partition. It is a transparent
Wagstaff-aligned conditional reference distribution, not an exact derivation
of the full discovery process.

The distribution is evaluated without Monte Carlo. For each interval form
\begin{equation}
G_b(z,t)
=
\prod_{p\in\mathcal P_b^-}(1+w_0(p)z)
\prod_{p\in\mathcal P_b^+}(1+w_0(p)zt).
\label{eq:weighted-subset-pgf}
\end{equation}
After normalization, the coefficient of \(z^{N_b}t^j\) is the exact
probability of selecting \(j\) higher-\(S\) exponents in interval \(b\).
Dynamic programming computes the required elementary symmetric sums, and
convolution of the six interval distributions gives the law of
\[
X_+=\sum_{b=1}^6|A_b\cap\mathcal P_b^+|.
\]

The historical value is \(X_+^{\rm obs}=32\). Under the primary coarse
weight, the exact conditional distribution has
\[
\mathbb E_0X_+=23.3279764307,
\qquad
\Pr_0(X_+\ge32)=0.00740204354948.
\]
For the centered, direction-agnostic convention
\[
p_{\rm two}
=
\Pr_0\!\left(
|X_+-\mathbb E_0X_+|
\ge |32-\mathbb E_0X_+|
\right),
\]
the result is
\begin{equation}
\boxed{p_{\rm two}=0.0115402626809.}
\label{eq:exact-randomization-p}
\end{equation}
The upper tail is the model-aligned directional contrast, but the two-sided
value is primary because the direction was not preregistered before the
historical pattern was examined.

The entire design was repeated with the more precise classical weight
\[
w_1(p)=\frac{\log(a(p)p)}p,
\qquad
a(p)=
\begin{cases}
2,&p\equiv3\pmod4,\\
6,&p\equiv1\pmod4.
\end{cases}
\]
After recomputing every interval partition, the observed split remains
\(32/15\). The refined distribution gives
\[
\mathbb E_1X_+=23.2346046013,
\quad
\Pr_1(X_+\ge32)=0.00675975252065,
\quad
p_{\rm two,1}=0.0111821205706.
\]
Thus the population-level association is not an artifact of replacing the
congruence-sensitive numerator by \(\log p\).

\subsection{Exact effective-modulus census and numerical evaluation}
\label{app:empirical}

The structural census contains two distinct operations: exact enumeration of
active pairs \((p,a)\), followed by numerical evaluation of the fully specified
local intensity \(\log F_a(p)\). They are audited separately.

\medskip
\noindent
\textbf{Exact active-pair enumeration.}
For every odd prime \(a<10^8\), the implementation computes
\(d=\operatorname{ord}_a(2)\) by exact modular exponentiation and visits the
progression
\[
p\equiv1\pmod d,
\qquad p>a,
\]
retaining the prime values with a least-prime-factor sieve and prime-index
array. This enumerates every pair
\[
13\le p<10^8,\qquad
a<p,\qquad
\operatorname{ord}_a(2)\mid p-1
\]
exactly once, without constructing or factoring \(A_p\).

The full run covers \(5{,}761{,}450\) prime exponents, performs
\(1{,}162{,}580{,}232\) progression checks, and records
\(130{,}844{,}312\) active pairs. Its count array agrees entry by entry with
an independently generated order-count census. For a fixed cutoff \(B\), the
exact diagnostic identity is
\begin{equation}
L_B(p)
=
\sum_{d\mid p-1}c_d\!\left(\min(B,p)\right),
\qquad
c_d(X)=\#\{a<X:a\ {\rm odd\ prime},\ \operatorname{ord}_a(2)=d\}.
\label{eq:appendix-cutoff-identity}
\end{equation}
The simpler expression \(\sum_{d\mid p-1}c_d(B)\) is valid only when
\(p\ge B\). The implementation always applies \(a<\min(B,p)\).

\medskip
\noindent
\textbf{Stable local-factor evaluation.}
For each exponent, the values in~\eqref{eq:wn} are evaluated in log space.
Terms more than \(44\) log units below the modal unnormalized weight are omitted;
the remaining weights are normalized by log-sum-exp. For \(a<1000\), the full
\(n\)-sum, including the gcd correction in~\eqref{eq:psing}, is tabulated on a
dense nonuniform grid in \(\sqrt{R_p}\) and linearly interpolated. A
deterministic sample of \(2{,}000\) direct evaluations gives maximum local
absolute error \(3.36\times10^{-6}\) and mean error
\(6.29\times10^{-8}\). Direct small-\(p\) test cases give a maximum cumulative \(H\)
discrepancy of \(1.76\times10^{-5}\); this is a finite-case check,
not a uniform population bound.

For \(a\ge1000\), the exponentially small gcd correction is retained in the
validation target and the main term is evaluated from the first four
cumulants of \(n-1\) under \(w_n(p)\). On the \(599\) admissible
\((p,a)\) pairs generated by \(600\) reproducible sampling attempts,
comparison with the complete \(n\)-sum gives maximum absolute
error \(6.12\times10^{-12}\), mean error approximately \(5.1\times10^{-14}\),
and 99th percentile approximately \(8.7\times10^{-13}\). On these checks, the errors are
several orders below the precision shown in the allocation tables.

The following deterministic test cases check the weight calculation, the exact active
set, and the resulting sector-averaged census calculation:
\begin{table}[H]
\centering
\small
\begin{tabular}{rrrrr}
\toprule
\(p\) & \(R_p\) & \(\sum_n n w_n\) & \(L(p)\) & \(H(p)\) \\
\midrule
101  & 10.9241711349 & 2.9039699848 & 5  & 1.4100215963 \\
1009 & 80.3093861742 & 8.2981717880 & 20 & 5.1061922209 \\
\bottomrule
\end{tabular}
\caption{Deterministic test cases for the unordered geometric-sector
calculation. The third column is a moment of the geometric sector measure,
not an estimate of the mean realized factorization length.}
\label{tab:factor-mass-fixtures}
\end{table}

The singleton formula is additionally compared with exhaustive enumeration of
every tuple in \(\Omega_{a,n}\) for \(a\in\{3,5,7,11\}\) and
\(2\le n\le5\). Every rational count agrees exactly.

\subsection{Finite-range singleton-geometry audit}
\label{app:singleton-audit}

The factor sample is regenerated from a dated Cunningham Project snapshot
\cite{CunninghamTables} whose SHA-256 digest is fixed in the archive. It
contains every completely factored composite \(M_p\) with prime \(p<1200\)
and at least three prime factors counted with multiplicity: \(143\) exponents
and \(833\) active pairs \((p,a)\). Listed multiplicities are retained, and
every factor product is
checked against \(2^p-1\).

The reference calculation does not approximate the number of singleton
coordinates by a binomial variable. For a multiplicity shape
\(\lambda=(\lambda_1,\ldots,\lambda_m)\vdash n\), let
\(s(\lambda)=\#\{j:\lambda_j=1\}\), and let \(m_t(\lambda)\) be the
multiplicity of the part size \(t\). The number of coordinate partitions with
this shape is
\[
\frac{n!}{\prod_j\lambda_j!\prod_t m_t(\lambda)!}.
\]
Writing the nonzero residues as powers of a generator of the cyclic group of
order \(h=a-1\), the number of distinct residue values satisfying the product
constraint is obtained by M\"obius inversion over set partitions. For a set
partition \(\pi\) of \(\{1,\ldots,m\}\), put
\[
c_B=\sum_{j\in B}\lambda_j,
\qquad
\mu(\pi)=\prod_{B\in\pi}(-1)^{|B|-1}(|B|-1)!.
\]
The injective count is
\[
J_h(\lambda)
=
\sum_{\pi}
\mu(\pi)\,\gcd\!\left(h,(c_B)_{B\in\pi}\right)h^{|\pi|-1}.
\]
Summing the product of these two counts over shapes with
\(s(\lambda)=s\) gives the exact distribution of
\(N_{\mathrm{sing}}(a,n)\). It sums to \(h^{n-1}\), and its first moment is
checked symbolically against
\(nP_{\mathrm{sing}}(a,n)\). Independent brute-force enumeration verifies the
complete distribution for \(a\in\{3,5,7\}\) and \(2\le n\le5\).

Within each exponent, CRT factors the uniform reference over distinct
active moduli. For the pooled statistic, we take the reference law to be
the product of these exponent-level laws and convolve the pair distributions.
For
\(T=\sum_{p,a}N_{\mathrm{sing}}(p,a)\), the pooled calculation gives
\[
T_{\mathrm{obs}}=2569,
\qquad
\mathbb E_\Omega T=2556.433352629154,
\qquad
\operatorname{sd}_\Omega(T)=22.991055146292.
\]
The centered two-sided reference value is
\[
\Pr_\Omega\!\left(
|T-\mathbb E_\Omega T|
\ge |T_{\mathrm{obs}}-\mathbb E_\Omega T|
\right)
=0.586877648165.
\]
None of the three disjoint subranges in Table~\ref{tab:singleton-diagnostic}
shows a detectable departure from the CRT reference. This audit compares
the singleton totals and secondary occupancy summaries; its reference
does not model complete-factorization selection, prime-factor size,
cross-sector weighting, modular aggregation, or transmission to primality.

As a secondary occupancy diagnostic, the archive also reconstructs the exact
distribution of the occupancy deficit \(d=n-m\). If
\(q_0(a,n)\) and \(q_1(a,n)\) denote its exact probabilities at \(d=0\) and
\(d=1\), respectively, the secondary offset model is
\[
\operatorname{logit}\Pr(d=1\mid d\in\{0,1\},a,n)
=
\alpha_n+\log\frac{q_1(a,n)}{q_0(a,n)}
+\gamma_{\mathrm{occ}}\log(a-1).
\]
The offset regression retains only pairs with \(3\le n\le6\),
observed \(d\in\{0,1\}\), and
\(q_0(a,n),q_1(a,n)>0\); pairs with structural-zero reference
probabilities are excluded. It uses \(564\) pairs from \(105\)
exponents and gives
\[
\widehat\gamma_{\mathrm{occ}}=-0.032564812626,
\qquad
\mathrm{CI}_{.95}^{\mathrm{cluster}}
= [-0.344057762228,0.195815197880].
\]
The all-collision variant uses the same factor-number range,
with outcome \(d\ge1\), and retains only pairs with
\(0<q_0(a,n)<1\); its offset is
\(\log[(1-q_0(a,n))/q_0(a,n)]\). It gives
\[
\widehat\gamma_{\mathrm{occ}}=0.030627369276,
\qquad
\mathrm{CI}_{.95}^{\mathrm{cluster}}
= [-0.260838334017,0.242655593528].
\]
Neither diagnostic detects a global occupancy tilt beyond the exact CRT
reference, and neither tests the subsequent transmission of the auxiliary
factor to Mersenne primality.

\subsection{Detailed versus reduced allocations}
\label{app:allocation-audit}

Table~\ref{tab:allocation-audit} reports two unanchored diagnostics and two
reference-anchored allocations. The detailed census diagnostic and the interval
reduced diagnostic use \(\log p+H(p)\) and
\(\log p+\theta_I\tau(p-1)\), respectively. The last two columns use the
relative \(S=1\) anchor. Every displayed total is conditioned
on the six observed interval counts; no structural coefficient uses the
higher-/lower-\(S\) labels.

\begin{table}[H]
\centering
\scriptsize
\begin{tabular}{lrrrrr}
\toprule
Interval & Census diag. & Reduced diag. & Anchored \(\theta_I\) & Anchored \(\theta\) & Observed high/total \\
\midrule
\([13,10^3)\)       & 5.403474 & 5.378090 & 5.424383 & 6.078077 & 7/10 \\
\([10^3,10^4)\)     & 4.515930 & 4.580234 & 4.691394 & 5.369126 & 5/8 \\
\([10^4,10^5)\)     & 3.547740 & 3.582089 & 3.722227 & 4.154373 & 4/6 \\
\([10^5,10^6)\)     & 3.108879 & 3.110838 & 3.285983 & 3.548235 & 4/5 \\
\([10^6,10^7)\)     & 3.258590 & 3.237232 & 3.485401 & 3.617595 & 3/5 \\
\([10^7,10^8)\)     & 8.776124 & 8.697482 & 9.544156 & 9.513498 & 9/13 \\
\midrule
Total                 & 28.610737 & 28.585964 & 30.153544 & 32.280904 & 32/47 \\
\bottomrule
\end{tabular}
\caption{Conditioned higher-\(S\) allocations. The first two columns are
unanchored diagnostics; the next two are reference-anchored model allocations.
The corresponding lower-\(S\) allocation is the interval total minus the entry
shown.}
\label{tab:allocation-audit}
\end{table}

For the interval-specific anchored model, the marginal factors
\(\beta_I\) are
\[
(.872923,.822918,.779858,.733424,.685702,.639282).
\]
Using the single global coefficient \(\theta=0.593255130346\), the interval factors
are
\[
(.651374,.647972,.648191,.644486,.643539,.642915),
\]
and the complete-population value is \(\beta=0.646960073626\). These factors
restore the Wagstaff marginal mass. They cancel from each within-interval
higher/lower share and are not fitted to the observed allocation.

The coefficient sensitivity is reported in
Table~\ref{tab:global-sensitivity}. In particular, weighting the projection
by \(\log p/p\) gives \(30.9040/16.0960\), the geometric per-modulus mean gives
\(33.2089/13.7911\), and the congruence-sensitive Wagstaff numerator,
with the same reference anchor and the refined-weight partitions,
gives \(31.5186/15.4814\). These alternatives show that the redistribution
amplitude is not uniquely identified.

As an additional baseline check, removing the universally active modulus
\(a=3\) reduces the global through-origin coefficient from \(0.593255\) to
\(0.583599\), while the anchored allocation changes only from
\(32.281/14.719\) to \(32.189/14.811\). Thus the global higher-\(S\)
displacement is not driven by the universal \(a=3\) contribution.

\paragraph{Aggregation-order sensitivity.}
The main benchmark marginalizes the geometric sector measure separately for
each active modulus before multiplying the resulting effective factors. As a
sensitivity check on this ordering, we also evaluated the alternative global-sector
construction
\[
\lambda_{\mathrm{shared}}(p)
=
\left[
\sum_{2\le n<R_p}
w_n(p)
\prod_{a\in\mathcal A_p}s_{a,n}
\right]^{-1},
\qquad
H_{\mathrm{shared}}(p)
=
\log\lambda_{\mathrm{shared}}(p).
\]
Here the same sector index \(n\) is retained across all active moduli before
marginalization. This alternative audits the order of sector marginalization
and modular aggregation.

Recomputing the complete population of \(5{,}761{,}450\) prime exponents with
the same active set, Wagstaff-balanced regions, and conditioned interval totals
gives an unanchored higher-/lower-\(S\) allocation of
\[
28.554607/18.445393,
\]
compared with \(28.610737/18.389263\) for the main census. The corresponding
global through-origin response coefficient is
\[
\theta_{\mathrm{shared}}=0.591976,
\]
compared with \(\theta=0.593255\). After applying the same \(S=1\) reference
anchoring and marginal Wagstaff calibration, the global allocation becomes
\[
32.268747/14.731253,
\]
compared with \(32.280904/14.719096\) under the main aggregation.

Thus retaining the sector index across moduli changes the unanchored
higher-\(S\) mass allocation by \(0.0561\) cases and the globally anchored
mass allocation by only \(0.0122\) cases. The two aggregation orders
therefore give effectively the same historical allocation.

\subsection{Sensitivity to the cross-sector measure}
\label{app:sector-prior-sensitivity}

The unordered geometric measure in~\eqref{eq:wn} allocates weight across
factor-number sectors \(n\); this cross-dimensional allocation is a modeling
choice rather than a counting theorem. To assess sensitivity to \(w_n\),
we tilt the same unordered sector weights by
\(n^{\gamma_{p,r}}\) while preserving their support and within-sector geometry.
With support endpoints \(n_{\min}(p)\) and \(n_{\max}(p)\), set
\(m_{p,r}:=\min\{n_{\max}(p),\max\{n_{\min}(p),r\mu_p\}\}\) and put
\begin{equation}
\mu_p=\sum_n n w_n(p),
\qquad
w_n^{(r)}(p)
=\frac{n^{\gamma_{p,r}}w_n(p)}
{\sum_m m^{\gamma_{p,r}}w_m(p)},
\qquad
\sum_n n w_n^{(r)}(p)=m_{p,r},
\label{eq:mean-shift-sector-prior}
\end{equation}
where \(r\in\{.75,.90,1,1.10,1.25\}\). For an interior \(m_{p,r}\), choose
finite \(\gamma_{p,r}\); at the lower or upper endpoint use the limiting
distribution \(\gamma_{p,r}\to-\infty\) or \(+\infty\). Endpoint projection
occurs only for \(14\) exponents at \(r=.75\), \(8\) at \(r=.90\), and
\(p=13,17\) for each upward shift.

The variants were recomputed over all \(5{,}761{,}450\) prime exponents
\(13\le p<10^8\) and all \(130{,}844{,}312\) active pairs \((p,a)\)
with \(a<p\), holding all other benchmark elements fixed. We report the
detailed unanchored allocation, the through-origin coefficient
\(\theta_r\) in \(H_r\simeq\theta_r\tau(p-1)\), and the allocation obtained
after the same \(S=1\) reference anchor and marginal calibration. The \(r=1\)
calculation independently reproduces the principal response to a maximum
pointwise difference below \(3.2\times10^{-8}\).

\begin{table}[H]
\centering
\small
\begin{tabular}{crrrr}
\toprule
\( r\mu_p \) & Unanchored high \(S\) & Unanchored low \(S\) & \(\theta_r\) & Anchored high/low \\
\midrule
\(.75\mu_p\) & 27.6990 & 19.3010 & .355198 & 29.6452/17.3548 \\
\(.90\mu_p\) & 28.4584 & 18.5416 & .571129 & 32.0679/14.9321 \\
\(\mu_p\) & 28.6107 & 18.3893 & .593255 & 32.2809/14.7191 \\
\(1.10\mu_p\) & 28.7464 & 18.2536 & .613626 & 32.4720/14.5280 \\
\(1.25\mu_p\) & 28.9311 & 18.0689 & .641422 & 32.7256/14.2744 \\
\bottomrule
\end{tabular}
\caption{Complete-population sensitivity to deliberate shifts of the mean
sector index. The first two allocation columns are unanchored diagnostics;
all allocations are conditioned on the same six interval totals, and the
anchored column uses the global coefficient \(\theta_r\).}
\label{tab:sector-prior-sensitivity}
\end{table}

Every variant remains displaced toward higher \(S\) relative to the balanced
\(23.5/23.5\) baseline. The \(25\%\) downward shift has the largest effect,
retaining \(4.20\) of the \(5.11\) mass-allocation displacement of the main
unanchored census. The spread in \(\theta_r\) and the anchored allocation
shows that the amplitude depends on the chosen sector measure.

\subsection{Extrapolation computation}
\label{app:prospective-computation}

The larger-interval population \(10^8\le p<10^9\) was generated independently by
a segmented prime sieve. It contains \(45{,}086{,}079\) prime exponents.
All extrapolated shares are evaluated using the fixed threshold
\(S(p)\ge S_{50}\) specified in Section~\ref{sec:prospective}, where
\[
S_{50}=1.071704.
\]
This rounded threshold gives the higher region coarse Wagstaff mass
\(0.499999645778\). No historical or future Mersenne-prime label enters
this boundary. Keeping the global \(\theta\) fixed, calibration over this
larger-interval population gives
\[
\beta=0.642691262463,
\qquad
q_+=0.738841185639,
\]
whereas continuation with the last measured interval coefficient
\(\theta_6 = 0.602695955689\) gives
\[
\beta=0.639057537448,
\qquad
q_+=0.741270134016.
\]
The coarse Wagstaff decade expectation \(5.9165945\) is therefore allocated by
the global model as \(4.3714237/1.5451708\). This absolute count retains the
uncertainty of the Wagstaff heuristic; the reported output of the fixed rule
is the higher-/lower-\(S\) share.

\subsection{Order-class decomposition and a post hoc 2-adic control}
\label{app:order-control}

The original higher-/lower-\(S\) comparison tests the full divisor
statistic \(S(p)\) under the design of
Appendix~\ref{app:wagstaff_balanced}. This additional post hoc control
examines its arithmetic interpretation and the reduction's fidelity under
a different partition. The benchmark, historical partition, and full-\(S\) extrapolation in
Section~\ref{sec:prospective} are unchanged.

\paragraph{Order classes and the exact census decomposition.}
Write \(p-1=2^v u\), where \(u\) is odd, and
\(d=\operatorname{ord}_a(2)=2^m t\), where \(t\) is odd. For an active prime
modulus \(a<p\), the activation condition \(d\mid p-1\) requires both
\(m\leq v\) and \(t\mid u\). We partition the active moduli into three
disjoint classes: pure powers of two (\(t=1\), \(m\geq1\)), pure odd orders
(\(m=0\), \(t>1\)), and mixed orders (\(m\geq1\), \(t>1\)). Writing
\(\mathcal A_{p,c}\) for the active moduli in class \(c\), define
\[
L_c(p)=|\mathcal A_{p,c}|,
\qquad
H_c(p)=\sum_{a\in\mathcal A_{p,c}}\log F_a(p).
\]
Then \(L(p)=\sum_c L_c(p)\) and \(H(p)=\sum_c H_c(p)\) identically, for
every exponent. The local factors, sector weights, and cutoff \(a<p\)
are unchanged; no Mersenne-prime label enters these responses.

\begin{table}[htbp]
\centering
\caption{Complete-population decomposition by multiplicative-order
class for \(13\leq p<10^8\) and \(a<p\). The last column attributes the
original unanchored higher-\(S\) \emph{mass-share displacement} within the complete
benchmark; it is not a causal attribution of observed primality.}
\label{tab:order-classes}
\begin{tabular}{lrrrr}
\toprule
Order class & Active pairs & Share of \(L\) & Share of \(\sum_p H\) & Displacement \\
\midrule
Pure power of two & \(11{,}619{,}521\) & \(8.88\%\) & \(12.54\%\) & \(+0.6306\) \\
Pure odd          & \(39{,}035{,}938\) & \(29.83\%\) & \(33.68\%\) & \(+1.6147\) \\
Mixed             & \(80{,}188{,}853\) & \(61.29\%\) & \(53.78\%\) & \(+2.8109\) \\
\midrule
Total             & \(130{,}844{,}312\) & \(100\%\) & \(100\%\) & \(+5.0563\) \\
\bottomrule
\end{tabular}
\end{table}

For clarity, the displacement column uses the actual finite-partition
Wagstaff mass share as reference: its higher-\(S\) total is \(23.5544\),
compared with \(28.6107\) for the detailed census. The rounded balanced
reference \(23.5/23.5\) used elsewhere in the manuscript is a different
summary.

For each prime exponent \(p\), set
\[
S_{\mathrm{odd}}(p)
=\frac{\log\tau((p-1)/2^{v_2(p-1)})}{\log\log p}.
\]
In the original six-interval Wagstaff-balanced design, without additional
\(v_2(p-1)\) stratification, the statistic
\(S_{2\mathrm{ad}}(p):=\log(v_2(p-1)+1)/\log\log p\) alone gives \(31/47\) in its higher
region (two-sided reference value \(0.0259\)); \(S_{\mathrm{odd}}\) alone gives
\(24/47\) (\(1.000\)). With the refined \(\log(a(p)p)/p\) weight, these
become \(30/47\) (\(0.0535\)) and \(24/47\) (\(0.8833\)), respectively.
Each post hoc comparison defines a different partition and does not
identify a causal 2-adic mechanism.

Under the original full-\(S\) partition, Table~\ref{tab:order-classes}
attributes \(87.53\%\) of the unanchored detailed-census displacement to
odd- and mixed-order classes. This model-based mass-share attribution and
the historical component comparisons concern different quantities; neither
is a causal decomposition of Mersenne primality.

\paragraph{The conditioned odd-part comparison.}
Within the six exponent intervals, stratify the complete prime-exponent
population by the exact value of \(v_2(p-1)\). There are 26 strata
containing at least one observed Mersenne-prime exponent. Preserve the
observed number of cases in every stratum. Within each stratum, order
prime exponents by decreasing \(S_{\mathrm{odd}}\) (breaking ties by
increasing \(p\)), and select the prefix closest to half of its total
Wagstaff weight \(\log p/p\) as its higher region. Within each stratum, the
ordering and half-weight boundary use only \(S_{\mathrm{odd}}\) and Wagstaff
weights; the stratum counts on which the comparison conditions are observed
Mersenne-prime counts. For each model, evaluate the same
without-replacement, item-weighted subset law and centered two-sided
reference convention as in Appendix~\ref{app:wagstaff_balanced}, conditioning on the count in
each of the 26 strata.

The corresponding unnormalized weights are \(\log p/p\),
\((\log p+H(p))/p\),
\([(1-\theta_I)\log p+\theta_I\tau(p-1)]/p\), and
\([(1-\theta)\log p+\theta\tau(p-1)]/p\), respectively. The marginal
factors \(\beta_I\) cancel within each conditioned stratum. The exact
subset means in Table~\ref{tab:odd-control} must not be confused with
totals obtained by multiplying each stratum's modeled mass share by
its observed count: without-replacement sampling makes these distinct.
The unrounded per-decade and per-stratum calculations are retained in
the reproducibility supplement; the decade-level contrasts are
supporting descriptions, not separately selected tests.

As a diagnostic, we rebuilt the within-stratum half-weight partition
using \(S\) while keeping the \((I,v_2)\) case counts fixed. This gives \(25/47\)
higher-region cases and a two-sided reference value
\(p_{\mathrm{two}}=0.7679\). Thus the absence of a separately
detectable association under this conditional design is not produced
solely by replacing \(S\) with \(S_{\mathrm{odd}}\).

\begin{table}[H]
\centering
\small
\caption{Post hoc 2-adic-stratified control. All rows use the same
higher-\(S_{\mathrm{odd}}\) regions and 25 observed cases among 47;
reported means and two-sided reference values are exact under each
specified weighted-subset law.}
\label{tab:odd-control}
\begin{tabular}{lrr}
\toprule
Conditional weight & Exact mean, higher region & Two-sided reference \\
\midrule
Wagstaff                                      & \(23.4245\) & \(0.6613\) \\
Unanchored detailed census                    & \(27.6417\) & \(0.4529\) \\
Interval-specific anchored reduction          & \(29.5543\) & \(0.1668\) \\
Globally anchored reduction                   & \(31.4765\) & \(0.0570\) \\
Unanchored detailed census, \(a<\min(1000,p)\)  & \(27.5027\) & \(0.4558\) \\
\bottomrule
\end{tabular}
\end{table}

\paragraph{Attribution within the new partition.}
For any stratum \(g\), let \(N_g\) be its fixed observed count,
\(b_p=\log p/p\), \(q_{0,g}=\sum_{g,+}b_p/\sum_g b_p\), and
\(D_g=\sum_g[b_p+\sum_c k_c(p)]\). The additive attribution of the
\emph{modeled mass-share displacement} is
\[
\Delta_{g,c}=\frac{N_g}{D_g}\left(
\sum_{p\in g,+}k_c(p)-q_{0,g}\sum_{p\in g}k_c(p)
\right).
\]
For the detailed census \(k_c(p)=H_c(p)/p\). With \(v_2(p-1)\) fixed,
the contributions summed over the 26 strata are \(+0.0020\) from pure
powers of two, \(+1.6749\) from pure odd orders, and \(+2.5607\) from mixed
orders, totaling \(+4.2376\) mass-share cases above the corresponding
Wagstaff shares.

To diagnose the reduced rule, we additionally project each response
component by the same through-origin criterion as
Equation~\eqref{eq:theta-projection}:
\[
\begin{aligned}
\theta_c
&=\frac{\sum_p\tau(p-1)H_c(p)}
{\sum_p\tau(p-1)^2},\\
(\theta_{\mathrm{pow2}},\theta_{\mathrm{odd}},\theta_{\mathrm{mix}})
&=(0.0521322,\,0.1865081,\,0.3546149),\\
\sum_c\theta_c&=0.5932551.
\end{aligned}
\]
Assigning the \emph{diagnostic projected components}
\(k_c(p)=\theta_c[\tau(p-1)-\log p]/p\) gives corresponding contributions
\(+0.7128\), \(+2.5500\), and \(+4.8485\), totaling \(+8.1113\) in
mass-share units. These are an algebraic attribution added for this
diagnostic, \emph{not} separately defined predictive components of the
original one-coefficient rule. In particular, since
\(\tau(p-1)=(v+1)\tau(u)\) at fixed \(v\), the projected pure-power-of-two
component inherits variation with the odd divisor count that is almost
absent in the original detailed \(H_{\mathrm{pow2}}(p)\) response.
This illustrates a loss of conditional fidelity of the projected and
anchored law; it does not demonstrate an error in
Equation~\eqref{eq:survival-an}.
Differences in these mass-share contributions cannot be directly
subtracted from differences between exact Fisher means.

This exploratory control finds no separately detectable odd-part association
under its conditional design. The order decomposition nevertheless shows
substantial odd- and mixed-order contributions to the \emph{constructed}
census. The globally anchored reduction has limited conditional fidelity
here, but the comparison does not isolate compression from anchoring or
establish a failure of the effective modular benchmark.

Supplement~S.3 reports a separate post hoc audit that retains the
detailed aligned even-order response without population regression
onto $\tau(p-1)$, contrasts raw and interval-centered transmission,
and conditions on $v_2(p-1)$. Its results do not alter the original
higher-/lower-$S$ inference or the frozen extrapolation.

\FloatBarrier

\section*{Declaration of Generative AI Use}

OpenAI ChatGPT (accessed July--November 2026)
assisted with English editing, generation and debugging of reproducibility
code, and computational cross-checks. The author reviewed and checked all
AI-assisted content, made all mathematical, methodological, and interpretive
decisions, and assumes full responsibility for the manuscript.

\section*{Data Availability Statement}
\footnotesize

GIMPS supplies the Mersenne-prime exponent data. Source code, the dated
Cunningham input snapshot, compact numerical references, and reproduction
instructions accompany this version of arXiv:2603.08994 as ancillary files.
The ancillary \path{README.md} gives the commands for the core calculations,
integrated controls, and the audits in Supplement~S.1, S.3, and S.4.
The \path{MANIFEST_SHA256.txt} file records distribution checksums, and
\path{VALIDATION.json} states the scope of the computational checks.

Large derived population arrays are regenerated by the included programs
in user-selected working directories. Further post hoc control details
and their code locations are indexed in
\path{archived_explorations/README.md}. These materials provide no
independent validation of the restricted transmission to primality.
The current article sources are \path{main.tex} and \path{supplement.tex}
at the submission root. In the code, the projection coefficient
\texttt{c} corresponds to \(\theta\) (or \(\theta_I\)), not the
Wagstaff constant \(C\).
\normalsize
\end{appendices}

\end{document}


\maketitle

This supplement documents benchmark anatomy and scope controls (S.1),
conditional information limits (S.2), the post hoc even-order and
\(2\)-adic audit (S.3), and residual-background sensitivity (S.4).
None derives transmission to primality or
establishes additional information beyond \(v_2(p-1)\).

\paragraph{Reproducibility.}
The ancillary \path{README.md} gives commands to reproduce the core
calculations and each computational audit. Source code, input data, and
compact numerical references are included; large working arrays are
regenerated. Further post hoc comparisons and their code locations are
indexed in \path{archived_explorations/README.md}.

\section{Benchmark anatomy and scope audits}
\label{app:identity-localization}

The fixed-channel density ratio \(1/(a-1)\) is independent of the channel.
The controls below retain the sector measure and distinguish
target-specific geometry, family selection, and the effects of modulus
size and sector weights. The homogeneous-localization audit is a separate
diagnostic and leaves the declared singleton response unchanged.

\paragraph{Exact local geometry for a general target.}
For an odd prime \(a\) and sector length \(n\ge2\), put \(m=a-1\).
For \(t\in\mathbb F_a^\times\),
define
\[
\Omega_{a,n}(t):=
\{(u_1,\ldots,u_n)\in(\mathbb F_a^\times)^n:\prod_i u_i=t\},
\qquad N_n(t):=\#\{x\in\mathbb F_a^\times:x^n=t\}.
\]
Every target has exactly \(m^{n-1}\) assignments. Under its uniform measure,
the singleton-coordinate probability is
\begin{equation}
P_{\mathrm{sing}}(a,n;t)=
\frac{(m-1)^{n-1}+(-1)^{n-1}(N_n(t)-1)}{m^{n-1}}.
\label{eq:general-target-singleton}
\end{equation}
Inclusion--exclusion on \(u_j=u_1\) gives intersection probabilities
\(m^{-j}\) until the final intersection, whose probability is
\(N_n(t)/m^{n-1}\), yielding Equation~\eqref{eq:general-target-singleton}.
For \(g=\gcd(n,m)\), \(N_n(t)=g\) if \(t^{m/g}=1\), and zero otherwise.
The case \(t=1\) recovers the main manuscript's formula.

Extend the same declared response by
\[
s_{a,n}(t)=P_{\mathrm{sing}}(a,n;t)
+\frac{1-P_{\mathrm{sing}}(a,n;t)}m,
\qquad F_{a,n}(t):=s_{a,n}(t)^{-1}.
\]
The target comparison is exact:
\[
P_{\mathrm{sing}}(a,n;1)-P_{\mathrm{sing}}(a,n;t)
=\frac{(-1)^{n-1}[g-N_n(t)]}{m^{n-1}}.
\]
If \(t\) is not an \(n\)-th power, \(F_{a,n}(1)>F_{a,n}(t)\) for even
\(n\), and \(F_{a,n}(1)<F_{a,n}(t)\) for odd \(n\). If it is an
\(n\)-th power, the responses coincide: multiplication of every coordinate
by a root identifies the fibers while preserving all coordinate collisions.
For fixed \(n\), the response differences are \(O_n(m^{1-n})\).
The sector-weighted difference is an alternating sum, so geometry
alone does not fix its sign. Exhaustive rational enumeration verifies
the formulas in \(158\) target-sector cases and \(1{,}444{,}100\) tuples.

\paragraph{Complementary homogeneous-channel localization: proof and finite-scale audit.}
Let $p\ge13$ be prime, let $a<p$ be an odd prime, and put
$t_a(p)=M_p\bmod a\in\mathbb F_a^\times$.
For the main manuscript's sector weights, define
\[
 w_*(p):=\max_{2\le n<R_p}w_n(p),\qquad
 W_a^{\mathrm{hom}}(p):=
 \max_{r\in\mathbb F_a^\times}
 \sum_{\substack{2\le n<R_p\\r^n=t_a(p)}}w_n(p).
\]
If $t_a(p)=1$, choosing $r=1$ gives $W_a^{\mathrm{hom}}(p)=1$.
For $t_a(p)\ne1$, put $d_a=\operatorname{ord}_a(t_a(p))$.
We first show that $d_a\ge4$. Order two would require
$t_a(p)=-1$, hence $2^p=0$ in $\mathbb F_a$, which is impossible.
If the target had order three, write $t=t_a(p)$ and $x=t+1=2^p$.
Then $t^2+t+1=0$, $x^2=t$, and $x^3=-1$, so $x$ has order six.
Since $\gcd(p,a-1)=1$, raising to the $p$-th power preserves
the order of $2$. This would force $\operatorname{ord}_a(2)=6$.
But $a\mid2^6-1=63$ would give $a=3$ or $7$, whose orders of
$2$ are two and three. This excludes order three as well.

For a fixed $r$ of order $\ell$, the compatible lengths, if any,
form one residue class modulo $\ell$, and $d_a\mid\ell$.
The sector weights are unimodal: on the consecutive positive support,
\[
 \frac{\widetilde w_{n+1}(p)}{\widetilde w_n(p)}
 =\frac{R_p-n-1}{n(n+1)}
   \left(1-\frac1{R_p-n}\right)^{n-1}
\]
decreases with $n$, since both positive factors decrease.
Each superlevel set $\{n:w_n(p)>u\}$ is therefore an interval.
A residue class modulo $\ell$ occupies at most $|J|/\ell+1$
points of an integer interval $J$. Integrating this count over
$0\le u<w_*(p)$ and using $\sum_nw_n(p)=1$ gives
\[
 \sum_{\substack{2\le n<R_p\\n\equiv c\!\!\pmod\ell}}w_n(p)
 \le \frac1\ell+w_*(p),\qquad
 W_a^{\mathrm{hom}}(p)\le\frac1{d_a}+w_*(p)
 \le\frac14+w_*(p).
\]
This is a deterministic statement for the declared sector measure.

If a consecutive sector window $I$ contains mass at least
$1-\varepsilon$, put $L_I=\max I-\min I$.
At most $1+\lfloor L_I/d_a\rfloor$ compatible lengths lie in $I$,
so
\[
 W_a^{\mathrm{hom}}(p)
 \le\left(1+\left\lfloor\frac{L_I}{d_a}\right\rfloor\right)
       w_*(p)+\varepsilon.
\]
Thus $d_a>L_I$ gives $W_a^{\mathrm{hom}}\le w_*+\varepsilon$,
and $2d_a>L_I$ gives $W_a^{\mathrm{hom}}\le2w_*+\varepsilon$.

The October 25 audit examines all complementary odd prime moduli
under $a<p/\log p$ for five exponents. For each exponent it takes
the smallest symmetric integer window about the sector mode whose
numerically computed mass is at least $1-10^{-8}$.
Table~\ref{tab:homogeneous-localization} counts the sufficient
order conditions above, not an enumeration of realized factorizations.

\begin{table}[htbp]
\centering\small
\setlength{\tabcolsep}{5pt}
\begin{tabular}{rrrr}
\toprule
$p$ & Complementary moduli & $d_a>L_I$ & $2d_a>L_I$ \\
\midrule
$1{,}000{,}003$ & $7{,}160$ & $7{,}049$ (98.45\%) & $7{,}103$ (99.20\%) \\
$10{,}000{,}019$ & $50{,}636$ & $50{,}434$ (99.60\%) & $50{,}533$ (99.80\%) \\
$50{,}000{,}017$ & $204{,}736$ & $204{,}425$ (99.85\%) & $204{,}558$ (99.91\%) \\
$90{,}000{,}049$ & $342{,}867$ & $342{,}506$ (99.89\%) & $342{,}660$ (99.94\%) \\
$100{,}000{,}007$ & $376{,}202$ & $375{,}791$ (99.89\%) & $375{,}980$ (99.94\%) \\
\bottomrule
\end{tabular}
\caption{Finite homogeneous-channel localization for $a<p/\log p$.
The last two columns guarantee at most one or two compatible
lengths, respectively, within the central window for any fixed $r$.
Counts and multiplicative orders are exact; window masses are
floating-point evaluations, without interval certification.}
\label{tab:homogeneous-localization}
\end{table}

For $p=90{,}000{,}049$, $w_*(p)\simeq0.0134495$;
approximately $99.895\%$ of the complementary moduli in this
diagnostic satisfy $W_a^{\mathrm{hom}}\le w_*+10^{-8}$.
These are finite observations, not an asymptotic proportion theorem.
The diagnostic cutoff $p/\log p$ is not the principal full census support
$a<p$; it has the order of the maximum sector length, whereas
the sector mode has order $\sqrt{R_p}$.
Fixing one homogeneous residue across sectors is distinct from
averaging the singleton response over every full fiber.
The bounds neither control the total response of the many
complementary moduli nor establish neutrality or a near-unit
residual multiplier in a primality law.

\paragraph{Selected-family generic response and modulus composition.}
The active family is
\(\mathcal A_p=\{a<p:a\text{ odd prime},\operatorname{ord}_a(2)\mid p-1\}\).
Let \(s^{\mathrm{gen}}_{a,n}\) denote the uniform target average of
\(s_{a,n}(t)\). Using the main manuscript's sector weights, define
\[
H=-\sum_{a\in\mathcal A_p}\log\!\left[\sum_n w_n(p)s_{a,n}(1)\right],
\quad
H_{\mathrm{ord}}=-\sum_{a\in\mathcal A_p}
\log\!\left[\sum_n w_n(p)s^{\mathrm{gen}}_{a,n}\right],
\quad
\Delta_{\mathrm{tgt}}=H-H_{\mathrm{ord}}.
\]
Uniform averaging over targets gives \(m^{-1}\sum_t N_n(t)=1\), since every
nonzero \(x\) contributes once to the pairs \((x,x^n)\). This gives the generic response and identity correction:
\[
s^{\mathrm{gen}}_{a,n}=\frac1m+\left(1-\frac1m\right)^n,
\qquad
s_{a,n}(1)-s^{\mathrm{gen}}_{a,n}
=(-1)^{n-1}\frac{(m-1)[\gcd(n,m)-1]}{m^n}.
\]

The control uses the same \(130{,}844{,}312\) active pairs.
In \(252\) deterministic prefix/small-modulus checks, the maximum
direct/interpolation difference is \(1.26\times10^{-5}\), not a uniform
bound. Both scores use the generic four-cumulant continuation for
\(a\ge1000\); its omitted target correction is not exactly zero and is
bounded below. The stored correction reaches \(0.843720\) at \(p=13\),
where only \(n=2\) has positive weight and the active moduli are \(3,5,7\).
Detailed reconstruction and approximation diagnostics remain archived.

For a fixed interval \(I\), higher region \(I_+\), and conditioned case
total \(N_I\), put \(B_I=\sum_I\log p/p\) and
\(\pi_{0,I}=\sum_{I_+}\log p/p\,/B_I\). The linear contrast is
\begin{equation}
D_I(G)=\sum_{p\in I_+}\frac{G(p)}p
-\pi_{0,I}\sum_{p\in I}\frac{G(p)}p,
\qquad
C(G)=\sum_I N_I\frac{D_I(G)}{B_I+\sum_{p\in I}H(p)/p}.
\label{eq:supplement-contrast}
\end{equation}
Write \(A_I(G):=N_Iq_I(G)\), where \(q_I(G)\) is the higher-region
mass share of \((\log p+G(p))/p\) when these weights are positive.
Each allocation uses its score's own
denominator; \(C\) uses the original \(H\) denominator for every component. Thus
\(C(H)=5.056296\), \(C(H_{\mathrm{ord}})=4.876617\), and
\(C(\Delta_{\mathrm{tgt}})=0.179679\): target averaging retains
\(96.45\%\) of the displacement from the finite-partition Wagstaff
allocation \(23.554441\). The target correction contributes \(0.096032\)
of \(D_I(H)=0.189318\) in the first interval; the aggregate percentage
is not a statement of uniform target independence. Universal identity
availability and the quantitative contribution of the product target are
distinct properties of the benchmark.

This counterfactual retains the selected modulus family. It neither
defines realized factorizations nor establishes primality relevance.

\begin{table}[htbp]
\centering\small
\setlength{\tabcolsep}{3.5pt}
\begin{tabular}{lrrrrrr}
\toprule
Interval & \(A_I(H)\) & \(A_I(H_{\mathrm{ord}})\) & \(D_I(H)\) & \(D_I(H_{\mathrm{ord}})\) & \(D_I(\Delta_{\mathrm{tgt}})\) & \(\theta_I^{\mathrm{gen}}\) \\
\midrule
\([13,10^3)\) & 5.403474 & 5.235342 & +0.189318 & +0.093286 & +0.096032 & 0.155139 \\
\([10^3,10^4)\) & 4.515930 & 4.515890 & +0.208571 & +0.208533 & +0.000038 & 0.234967 \\
\([10^4,10^5)\) & 3.547740 & 3.547740 & +0.331712 & +0.331712 & \(+2.79\!\times\!10^{-10}\) & 0.308550 \\
\([10^5,10^6)\) & 3.108879 & 3.108879 & +0.492676 & +0.492676 & \(0^\dagger\) & 0.390900 \\
\([10^6,10^7)\) & 3.258590 & 3.258590 & +0.686403 & +0.686403 & \(0^\dagger\) & 0.490918 \\
\([10^7,10^8)\) & 8.776124 & 8.776124 & +0.888480 & +0.888480 & \(0^\dagger\) & 0.602696 \\
\bottomrule
\end{tabular}
\caption{Target averaging within the same selected family. Allocations
\(A_I\) use each score's own denominator; \(D_I\) add on a common linear
scale. \(\dagger\) The target-specific term is numerically unresolved
at the displayed precision. The continuation omits it; its exact magnitude
is bounded separately below, rather than asserted to be zero.}
\label{tab:selected-target-average}
\end{table}

The global compression coefficient changes from \(0.593255130346\) to
\(0.593255118862\); in the smallest interval it changes from
\(0.161363\) to \(0.155139\).

To distinguish modulus count from composition, fit without outcome labels
\[
\eta_I=\frac{\sum_{p\in I}L(p)H_{\mathrm{ord}}(p)}{\sum_{p\in I}L(p)^2},
\qquad H_L=\eta_I L,\qquad
\Delta_{\mathrm{comp}}:=H_{\mathrm{ord}}-H_L.
\]
Then \(H=H_L+\Delta_{\mathrm{comp}}+\Delta_{\mathrm{tgt}}\), and
\(D_I(H_{\mathrm{ord}})=D_I(H_L)+D_I(\Delta_{\mathrm{comp}})\).

With the original score's denominators, count contributes \(5.201400\),
composition \(-0.324783\), and the target correction \(0.179679\),
summing to \(5.056296\).
The modulus-count component therefore accounts for the overall direction of the benchmark displacement; composition moderates it, although its contrast is positive in the last two intervals.
The full interval decomposition is archived. The count-only normalized allocation is \(29.116146\), but allocation
differences are not additive. The composition residual can yield negative individual weights and is
therefore used only in the additive decomposition.

\paragraph{A matched extension over auxiliary moduli.}
Let \(\mathcal B_p=\{a<p:a\text{ odd prime}\}\), and put
\(t_a(p)=M_p\bmod a\). This target is nonzero because prime divisors of
\(M_p\) are at least \(2p+1\). Retaining the same sector weights and
ordering, define
\[
F_a(p;t)=\left[\sum_n w_n(p)s_{a,n}(t)\right]^{-1},
\qquad
H_{\mathrm{all}}(p)=\sum_{a\in\mathcal B_p}\log F_a(p;t_a(p)).
\]
The selected-family component is \(H\), since \(t_a(p)=1\) exactly
for \(a\in\mathcal A_p\). For the matched complement,
\begin{equation}
H_{\mathrm{compl}}(p)=
\sum_{a\in\mathcal B_p\setminus\mathcal A_p}\log F_a(p;t_a(p)),
\qquad H_{\mathrm{all}}=H+H_{\mathrm{compl}}.
\label{eq:completion-score-decomposition}
\end{equation}
Here \(H_{\mathrm{compl}}\) is the score of the complementary moduli,
whereas \(\Delta_{\mathrm{comp}}\) is the residual due to modulus composition.
The declared logarithmic indices are additive; neither component is
assigned a separate primality contribution.

\paragraph{Reproduced extension and numerical checks.}
The unanchored selected-family allocation is \(28.610737\).
The matched all-moduli extension gives approximately \(23.0413\),
against the finite-partition Wagstaff allocation \(23.5544\).
For \(a<\min(1000,p)\), the corresponding allocations are
\(28.470979\) and approximately \(23.1382\). These are conditioned
mass-share allocations, not exact weighted-subset means.

Realized targets are evaluated directly for \(p<20{,}000\).
Beyond this prefix, the retained response is averaged uniformly over
targets using \(s^{\mathrm{gen}}_{a,n}\) above.
This is not an average of inverse responses. The bound
\(\lvert N_n(t)-1\rvert\le n-1\) controls the omitted target term.
With \(s^{\mathrm{gen}}_a\ge1/m\), its relative magnitude is at most
\((m-1)\sum_n w_n(p)(n-1)m^{1-n}\).
If this bound is \(\epsilon_a<1\), the absolute logarithmic error is
at most \(-\log(1-\epsilon_a)\).
The sector measure shifts monotonically toward larger \(n\) as \(R_p\)
increases, and \((n-1)m^{1-n}\) is nonincreasing in \(n\). Evaluating the
finite-modulus bound at the minimum radius over the studied population
\(20{,}000\le p<10^8\) gives a uniform omitted-score bound of \(1.094\times10^{-7}\), separately from
interpolation and cumulant errors.

The generic sum uses the supplied numerical grid, four cumulants, and a
fourth-order \(1/(a-1)\) expansion for \(a\ge10^5\).
Its finite direct/interpolation and cumulant checks do not provide
uniform roundoff bounds; complete diagnostics remain archived.

\paragraph{Contrasts on a common scale and a size control.}
Use the common-scale contrast \(D_I\) defined in
Equation~\eqref{eq:supplement-contrast}. Equation~\eqref{eq:completion-score-decomposition} gives
\(D_I(H_{\mathrm{all}})=D_I(H)+D_I(H_{\mathrm{compl}})\)
exactly. To separate reference geometry from realized targets, define
\[
H_{\mathrm{gen}}(p):=
-\sum_{a\in\mathcal B_p}\log\!\left[\sum_nw_n(p)s^{\mathrm{gen}}_{a,n}\right],
\qquad
\Delta_{\mathrm{all,tgt}}:=H_{\mathrm{all}}-H_{\mathrm{gen}}.
\]
This control retains the cutoff and the reference sector radius, including
its \(p\bmod4\) dependence, but uses neither realized targets nor the
identity-active family. With \(b(p)=\log p/p\),
\(D_I(G)=B_I\operatorname{Cov}_{b,I}
(\mathbf1_{I_+},G/\log p)\), where the covariance uses the normalized
weights \(b(p)/B_I\) on interval \(I\). A descriptor depending only on reference
size and residue class need not have zero contrast across these regions.

\begin{table}[htbp]
\centering\small
\setlength{\tabcolsep}{4pt}
\begin{tabular}{lrrrrr}
\toprule
Interval & \(D_I(H)\) & \(D_I(H_{\mathrm{all}})\) & \(D_I(H_{\mathrm{compl}})\) & \(D_I(H_{\mathrm{gen}})\) & \(D_I(\Delta_{\mathrm{all,tgt}})\) \\
\midrule
\([13,10^3)\) & +0.189318 & -0.148503 & -0.337820 & -0.238458 & +0.089956 \\
\([10^3,10^4)\) & +0.208571 & -0.045268 & -0.253839 & -0.045306 & \(+3.87\!\times\!10^{-5}\) \\
\([10^4,10^5)\) & +0.331712 & -0.209811 & -0.541524 & -0.209811 & \(+2.77\!\times\!10^{-10}\) \\
\([10^5,10^6)\) & +0.492676 & -0.358174 & -0.850850 & -0.358174 & \(0^\dagger\) \\
\([10^6,10^7)\) & +0.686403 & -0.105001 & -0.791404 & -0.105001 & \(0^\dagger\) \\
\([10^7,10^8)\) & +0.888480 & -1.864960 & -2.753440 & -1.864960 & \(0^\dagger\) \\
\bottomrule
\end{tabular}
\caption{Matched higher-\(S\) contrasts and generic size control.
Scores add on this common scale; normalized allocations do not.
\(\dagger\) The target-specific term is numerically unresolved at the
displayed precision and omitted by the generic continuation. Its exact
contrast is bounded by \(1.0\times10^{-8}\) in each of these intervals,
separately from other numerical errors.}
\label{tab:completion-contrasts}
\end{table}

The generic size/sector descriptor gives negative contrasts in all
six intervals and a conditioned allocation of \(22.903389\), compared
with \(23.041339\) for the realized-target extension.
Direct re-evaluation of all \(2{,}257\) exponents below \(20{,}000\)
reproduces both scores. Beyond this prefix, the omitted-score bound gives
\[
|D_I(\Delta_{\mathrm{all,tgt}}^{\mathrm{omitted}})|
\le1.094\times10^{-7}
\sum_{\substack{p\in I\\p\ge20{,}000}}
\frac{|\mathbf1_{I_+}(p)-\pi_{0,I}|}{p}.
\]
In \([10^7,10^8)\), the bound is \(7.31\times10^{-9}\), whereas
\(D_I(H_{\mathrm{all}})=-1.864960\).
The selected-family and all-moduli denominators are approximately
\(5.0745\) and \(201.1936\), respectively. Thus a near-balanced
allocation does not establish \(S\)-neutrality. This extension is a
scope control for the declared benchmark, not a derived complete sieve
or a determination of which response should enter a primality law.

For an actual factorization, every auxiliary product congruence follows
from the product identity. Extending the declared response to all these
moduli does not derive additional independent primality filters or
establish cancellation of the selected arithmetic contribution.

\paragraph{Reproduced alignment controls.}
For \(a<\min(1000,p)\), transfer a nearby donor's active family to
each recipient, retaining the recipient's sector response and cutoff.
Donors and recipients lie in the same exponent interval. Within
prime-rank blocks of sizes \(256,1024,4096\), use \(100\) fixed affine
bijections each (PCG64 seed \(20261009\)); partial blocks retain their
actual lengths. The code and numerical references are included in the ancillary files.

\begin{table}[htbp]
\centering\small
\begin{tabular}{lr}
\toprule
Control & Original \(S\) partition \\
\midrule
Wagstaff finite-partition reference & 23.554441 \\
Aligned selected family, \(a<1000\) & 28.470979 \\
Borrowed sets, block 256; mean & 23.577561 \\
Borrowed sets, block 1024; mean & 23.566974 \\
Borrowed sets, block 4096; mean & 23.572187 \\
\bottomrule
\end{tabular}
\caption{Borrowed-family alignment control on the conditioned 47-case
scale. Entries average 100 fixed affine permutations per block size;
the recipient's sector response and cutoff are retained. This is a
benchmark diagnostic, not a test of transmission to primality.}
\label{tab:preliminary-localization}
\end{table}

Borrowing active sets removes most of the displacement in this
design. The contrast depends on pairing an exponent with its own family,
beyond the retained sector response and cutoff; this does not establish
transmission to primality.

\FloatBarrier

\section{Illustrative conditional comparisons and information limits}
\label{app:prospective-evidence}

On any specified, completely ascertained eligible population \(\mathcal P\),
compare the
weights \(b_p=\log p/p\) with the principal redistribution weights
\(v_p=[(1-\theta)\log p+\theta\tau(p-1)]/p\).
The illustration compares the principal \(k=1\) rule with Wagstaff. Condition on
\(N\) distinct cases. For non-negative weights \(w\), let
\(e_N(w)=\sum_{|T|=N}\prod_{p\in T}w_p\), with \(e_0=1\), and require
\(e_N(w)>0\). The conditional
reference law and full-set likelihood ratio are
\begin{equation}
P_w(T\mid N)=\frac{\prod_{p\in T}w_p}{e_N(w)},
\qquad
\mathrm{LR}(T)=\frac{e_N(b)}{e_N(v)}\prod_{p\in T}\frac{v_p}{b_p}.
\label{eq:prospective-exact-lr}
\end{equation}
Common calibration scalars cancel exactly. If only the higher-region count
\(K\) is retained, use the ratio of the two count probabilities
\[
P_w(K=k\mid N)=
\frac{e_k(w_+)e_{N-k}(w_-)}{e_N(w)}.
\]
The full-set and count likelihood ratios retain different information.
Separately conditioned interval ratios multiply. A prospective application
requires an eligible population and documented search exposure.

For the complete population \(10^8\le p<10^9\), the descriptive
mass-share approximation is
\[
\mathrm{LR}_{\mathrm{bin}}(K,N)=
\left(\frac q{\pi_0}\right)^K
\left(\frac{1-q}{1-\pi_0}\right)^{N-K},
\qquad q=0.738841,\quad\pi_0=0.49999965.
\]
This approximation does not replace the exact without-replacement law.
For another eligible cohort, recompute \(q\) and \(\pi_0\) from frozen
weights and its partition. For six cases, only \(6/6\) gives a ratio
above ten (\(10.4109\)); its probability under the stated model is
\(q^6\approx0.1627\). Both \(0/6\) and \(1/6\) give ratios below
\(0.1\). For one case, the ratios are approximately \(1.48\) and
\(0.52\); conditioned on \(N=0\), the ratio is one.
These comparisons do not validate the amplitude or exclude weaker
modulations. The complete illustrative table remains archived.

The supplied \path{repro_v1.9.9/v10/audit_prospective_evidence.py} checks the exact
formulas by subset enumeration and reproduces the table.
An independent prospective comparison requires a frozen rule, a fully
ascertained eligible cohort, and documented search exposure. No such
observation or validation is claimed here.

\FloatBarrier

\section{Post hoc direct transmission of aligned even orders}
\label{sec:even-direct}

This exploratory audit leaves the principal rule, its $S=1$ reference,
and its fixed extrapolation unchanged. It instead examines
whether part of the historical pattern can be reproduced without projecting
the detailed modular score onto $\tau(p-1)$.
Define the aligned even-order subfamily and its uncompressed score by
\[
\begin{gathered}
\mathcal A_{p,\mathrm{even}}:=
 \{a<p:\ a\text{ an odd prime},\ 
 2\mid\operatorname{ord}_a(2),\ 
 \operatorname{ord}_a(2)\mid p-1\},\\[3pt]
H_{\mathrm{even}}(p):=
 \sum_{a\in\mathcal A_{p,\mathrm{even}}}\log F_a(p).
\end{gathered}
\]
The same local singleton response $F_a(p)$, sector measure and cutoff $a<p$
as in the main paper are retained. Even order includes mixed orders
$2^j d$ with $d>1$ odd. The condition
$\operatorname{ord}_a(2)\mid p-1$ is essential: this is not the
family of \emph{all} even-order auxiliary moduli.

For the six original exponent intervals $I$, let
\[
b(p)=\frac{\log p}{p},\qquad
c_I=\frac{\sum_{p\in I}H_{\mathrm{even}}(p)/p}
          {\sum_{p\in I}\log p/p},\qquad
w_{\mathrm{raw}}(p)=\frac{\log p+H_{\mathrm{even}}(p)}p.
\]
\[
w_{\mathrm{dir},I}(p)=
\frac{\log p+H_{\mathrm{even}}(p)-c_I\log p}{p}.
\]
The latter rule conserves the Wagstaff mass
$\sum_{p\in I}w_{\mathrm{dir},I}(p)=\sum_{p\in I}b(p)$ exactly.
The observed minimum multiplier $w_{\mathrm{dir},I}(p)/b(p)$ is
$0.243335$ on the historical prime-exponent population. This is a
population-mean centering, not the original $S=1$ reference anchor.
Neither its positivity nor predictive adequacy is established outside
the finite audited range.

Write \(S_1(p)=S(p)\) and define the purely \(2\)-adic partition by
\[
S_{1,\mathrm{2ad}}(p)=\frac{\log(1+v_2(p-1))}{\log\log p}.
\]
The main manuscript denotes the latter by \(S_{2\mathrm{ad}}\).

\begin{table}[htbp]
\centering\small
\caption{Retrospective higher-region mass allocations, conditional on the
six interval counts, on the 47-case scale. They are not exact
without-replacement means.}
\label{tab:even-direct-shifts}
\begin{tabular}{lrrrr}
\hline
Partition & Observed & Wagstaff & Raw & Centered\\
\hline
$S_1$ & 32 & 23.5544 & 27.6122 & 30.1115\\
$S_{1,\mathrm{2ad}}$ & 31 & 23.4266 & 27.5390 & 29.9788\\
\hline
\end{tabular}
\end{table}

The same $p-1$-selected response is used for both partitions, without
fitting a separate coefficient. Under the exact within-interval Fisher
weighted-subset law, the mean higher-$S_1$ counts are $23.32798$
(Wagstaff), $27.34070$ (raw), and $29.82817$ (centered).
The two-sided historical reference probability for the latter model at
32 observed cases is $0.53364$. This retrospective compatibility is
not independent predictive validation.

\paragraph{Pure and mixed even orders.}
The even-order active-pair census contains $91{,}808{,}374$ pairs
$(p,a)$, separated into $11{,}619{,}521$ pairs with pure
power-of-two order and $80{,}188{,}853$ with mixed order $2^jd$
($d>1$ odd). With componentwise centering, the modeled higher-$S_1$
mass displacements above Wagstaff are $+1.0511$ and $+5.5060$ cases,
respectively. Thus about $84\%$ of the \emph{constructed benchmark's}
direct-centering displacement arises from mixed even orders. This
percentage is not a causal attribution of the historical primality cases.

\paragraph{Conditioning on the $2$-adic valuation.}
In the 26 populated $(I,v_2(p-1))$ strata of the main manuscript's
Appendix A.14, preserve the observed case count in each stratum and
rebuild the higher/lower partition using
$S_{\mathrm{odd}}(p)=\log\tau((p-1)/2^{v_2(p-1)})/\log\log p$.
The observed count is $25$ higher-region cases out of $47$.
The exact Fisher means are $23.42447$ under Wagstaff,
$26.39737$ under the raw even-order rule, and $28.08035$ under the
centered rule. The corresponding log-likelihood ratios for the
\emph{observed individual exponents}, conditioned on the 26 stratum
counts and compared with Wagstaff, are $+0.9608$ (raw) and $-0.1623$
(centered). They are not $p$-values. This sample does not identify a
robust residual primality association beyond $v_2(p-1)$.

\paragraph{Additional valuation-threshold audit.}
The exploratory event $v_2(p-1)\ge3$ occurs in $19$ of the $47$
historical cases. The simple $47/4=11.75$ reference assumes
an idealized quarter-mass for the class $p\equiv1\pmod8$; it is
\emph{not} the exact Appendix A.8 reference. Conditioning on the
original six case counts and using the item-weighted
without-replacement law gives an exact mean $11.360997$, upper tail
$0.009183$, and centered two-sided probability $0.010499$ under
$b(p)=\log p/p$. Replacing the baseline by Wagstaff's refined
$\log(a(p)p)/p$ gives mean $11.905152$, upper tail $0.015211$,
and centered two-sided probability $0.018275$.
The threshold and the order restriction were examined post hoc; these
correlated reference calculations are not independent confirmatory tests.

\paragraph{Scope and sensitivity.}
Restricting the even-order family to $a<\min(1000,p)$ gives a centered
higher-$S_1$ allocation of $29.7321$ instead of $30.1115$.
Varying the direct intensity changes full-set likelihoods; no
intensity is selected by optimization over historical labels.
Shifted-partition comparisons and the complete original audit remain
available in \path{archived_explorations/even_order_full_audit.md}.

\paragraph{Interpretation and provenance.}
The local response remains a declared benchmark and its transmission
to Mersenne-primality weights is heuristic. Neither the 47 historical
cases nor the conditioned comparisons demonstrate that the detailed
modular score carries detectable primality-related information beyond
the $2$-adic valuation. The numerical large-modulus approximations are
audited on finite fixtures, not certified by a population-uniform error
bound; the alternative rule lacks a full sector-prior and prospective
sensitivity analysis. This audit is dated \textbf{19 October 2026},
and does not modify the principal transmission rule or extrapolation.

\FloatBarrier
\clearpage

\section{Post hoc residual-background sensitivity}
\label{sec:residual-background-sensitivity}

This post hoc analysis perturbs the principal rule by a positive residual
proxy. Eleven shifted responses serve only to construct this proxy,
without being proposed as alternative primality models. For
\(k=1,\ldots,12\), put
\[
 r_k(p)=1-\vartheta_k+\vartheta_k\frac{\tau(p-k)}{\log p},
 \qquad
 Z(p)=\frac1{11}\sum_{k=2}^{12}\log r_k(p),
\]
and define
\[
 G_\eta(p)=\exp(\eta Z(p)),\qquad
 w_\eta(p)=\frac{\log p}{p}\,r_1(p)G_\eta(p).
\]
Here \(\vartheta_1=\theta=0.593255130346\) is the principal full-census
\(a<p\) coefficient. The fixed \(\vartheta_k\), \(k=2,\ldots,12\), come
from label-independent shifted-family benchmarks with
\(a<\min(1000,p)\), not the principal cutoff. Their exact values are
supplied in \path{analyses_2026-10-17/coeficientes_Hk_fuente.csv}; its
\(k=1\) coefficient is not used. All \(r_k\) are positive on the audited
population. The strength \(\eta\) is chosen, not estimated from primality.

The original higher-/lower-\(S\) partition, six interval counts totaling
47 cases, and population of \(5{,}761{,}450\) prime exponents
\(13\le p<10^8\) are fixed. Each allocation sums the interval case count
times its higher-region weight share. Positive interval normalization
factors cancel. These mass allocations are not exact
without-replacement means. At \(\eta=0\), \(G_0=1\) recovers the principal
rule; this identity asserts no neutrality of actual omitted arithmetic.

\begin{table}[htbp]
\centering\small
\begin{tabular}{lcc}
\toprule
\(\eta\) & Higher-\(S\) allocation & Lower-\(S\) allocation \\
\midrule
\(-2\) & 33.093 & 13.907 \\
\(-1\) & 32.691 & 14.309 \\
\(0\) (principal rule) & 32.281 & 14.719 \\
\(1\) & 31.863 & 15.137 \\
\(2\) & 31.438 & 15.562 \\
\(11\) (strong stress) & 27.310 & 19.690 \\
\midrule
Observed & 32 & 15 \\
\bottomrule
\end{tabular}
\caption{Post hoc residual-proxy sensitivity on the conditioned 47-case
scale. Calculation dated October 17, 2026.}
\label{tab:selected-reference-sensitivity}
\end{table}

The tested window \(\eta=-2,-1,0,1,2\) and its allocation range
31.438--33.093 are illustrative scenarios, not confidence bounds or
estimates of the unresolved arithmetic contribution. Negative \(\eta\)
reverses this proxy without an assigned arithmetic interpretation.
The stress \(\eta=11\) gives \(G_{11}=\prod_{k=2}^{12}r_k\), a deliberately
strong product, not a derived joint transmission law. Neither the proxy
nor this window bounds total model uncertainty or independently validates
the principal rule. The code, coefficient inputs, and reference results
are supplied in \path{analyses_2026-10-17/}.